\documentclass[11pt]{article}
\usepackage{mathrsfs}
\usepackage{amsfonts}
\usepackage{natbib}
\usepackage{amsmath}
\usepackage{amsthm}
\usepackage{amssymb}
\usepackage{color}
\usepackage{lscape}
\usepackage{rotating}
\usepackage{caption2}
\usepackage{subfigure}
\usepackage{graphicx}

\usepackage{amssymb}

\def\be{\begin{equation}}
\def\ee{\end{equation}}
\def\bea{\begin{eqnarray}}
\def\eea{\end{eqnarray}}
\def\bd{\begin{displaymath}}
\def\ed{\end{displaymath}}
\def\bda{\begin{eqnarray*}}
\def\eda{\end{eqnarray*}}
\def\bsm{\begin{small}}
\def\esm{\end{small}}

\def\t0{\theta_0}

\def\ha1{\hat \beta_1}

\def\bnt{\begin{enumerate}}
\def\ent{\end{enumerate}}
\def\T{{ \mathrm{\scriptscriptstyle T} }}

\def\bsc{\begin{scriptsize}}
\def\esc{\end{scriptsize}}

\def\half{\txs{1\over 2}}

\newtheorem{theorem}{Theorem}

\newtheorem{ft}{Fact}[section]
\theoremstyle{definition}
\newtheorem{as}{Assumption}

\makeatletter
\newcommand{\figcaption}{\def\@captype{figure}\caption}
\newcommand{\tabcaption}{\def\@captype{table}\caption}
\makeatother

\newcommand{\cov}{{\rm Cov}}

\newcommand{\etal}{\mbox{\sl et al.\;}}

\newcommand{\var}{\mbox{Var}}

\def\ga{\gamma}

\newcommand{\bX}{{\mathbf X}}

\newcommand{\bc}{{\mathbf c}}

\def\hth{\hat\theta}

\renewcommand{\theequation}{\arabic{equation}}

\def\T{{ \mathrm{\scriptscriptstyle T} }}

\def\ab{\allowbreak}

\def\pr{{\mathrm{pr}}}
\def\a{\alpha}

\def\A{{\rm A}}
\def\ab{\allowbreak}
\def\as{^*}
\def\asas{^{**}}
\def\bc{^{{\rm bc}}}
\def\bcc{^{{\rm bcc}}}
\def\be{\beta}
\def\bX{{\bar X}}
\def\bigmi{\;\big|\;}

\def\biggmi{\;\bigg|\;}
\def\Biggmi{\;\,\Bigg|\;\,}
\def\cE{{\cal E}}
\def\cI{{\cal I}}

\def\cX{{\cal X}}
\def\cf{^{{\rm cf}}}
\def\cov{{\rm cov}}

\def\etal{et al.~}

\def\ga{\gamma}
\def\hA{{\hat A}}
\def\half{^{1/2}}
\def\hbe{{\hat\be}}

\def\hF{{\widehat F}}
\def\hG{{\widehat G}}

\def\heta{{\hat\eta}}

\def\hp{{\hat p}}
\def\hQ{{\widehat Q}}

\def\hth{{\hat\th}}
\def\hsi{{\hat\si}}
\def\hx{{\hat x}}
\def\hxi{{\hat\xi}}

\def\IR{{\rm I\!R}}
\def\ji{_{ji}}

\def\mhf{^{-1/2}}
\def\mo{^{-1}}
\def\mt{^{-2}}
\def\mth{^{-3}}
\def\mi{\,|\,}
\def\ooB{{1\over B}}
\def\oon{{1\over n}}
\def\osx{{\textstyle{1\over6}}}
\def\otf{{\textstyle{1\over24}}}
\def\part{\partial}
\def\ra{\to}
\def\rai{\ra\infty}

\def\si{\sigma}

\def\sumion{\sum_{i=1}^n\,}

\def\sumjop{\sum_{j=1}^p\,}

\def\tth{{\tilde\th}}
\def\tF{{\widetilde F}}
\def\th{\theta}
\def\thf{{\textstyle{1\over2}}}
\def\var{{\rm var}}

\begin{document}

\title{{\bf Double-bootstrap methods that use a single double-bootstrap simulation}}

\author{
Jinyuan Chang \qquad Peter Hall\\Department of Mathematics and Statistics\\ The University of Melbourne, VIC, 3010, Australia\\ jinyuan.chang@unimelb.edu.au~~~~~halpstat@ms.unimelb.edu.au}

\date{}
\maketitle

\begin{abstract}
We show that, when the double bootstrap is used to improve performance of bootstrap methods for bias correction, techniques based on using a single double-bootstrap sample for each single-bootstrap sample can be particularly effective.  In particular, they produce third-order accuracy for much less computational expense than is required by conventional double-bootstrap methods.  However, this improved level of performance is not available for the single double-bootstrap methods that have been suggested to construct confidence intervals or distribution estimators.
\end{abstract}

\noindent{\bf Keywords:} Bias correction; Bias estimation; Confidence intervals; Distribution estimation; Edgeworth expansion; Second-order correctness; Third-order correctness.

\setcounter{equation}{0}
\section{Introduction}

Double-bootstrap methods that use a single simulation at the second bootstrap level have been studied in at least one context for more than a decade.  An early contribution was made by \cite{White_2000}, although in the setting of diagnosing the overuse of a dataset, rather than speeding up Monte Carlo simulation for general applications of the bootstrap. \cite{DavidsonMackinnon_2001,DavidsonMackinnon_2002}, and the same authors in a number of subsequent papers accessible via \cite{Mackinnon_2006} and \cite{DavidsonMackinnon_2007}, introduced the concept independently and explored its applications. \cite{Giacominietal_2013} christened the technique the warp-speed double-bootstrap method, nomenclature that we shall use here, too. \cite{Giacominietal_2013} demonstrated that this approach is asymptotically consistent.  All this work is for the case of distribution estimation and its application to constructing confidence intervals and hypothesis tests.

In statistics the conventional double bootstrap is used in two main classes of problems: (i)~To improve the effectiveness of bias correction, and (ii)~to improve the coverage accuracy of confidence intervals.  In problem~(i), an application of the double bootstrap reduces the order of magnitude of bias by the factor~$O(n\mo)$, and in problem~(ii) it reduces coverage error by the factor $O(n\mhf)$ for one-sided confidence intervals, and $O(n\mo)$ for two-sided intervals.  In the setting of problem~(i), it is not clear whether there exists a version of warp-speed methodology for bias correction, and whether, should it exist, it successfully reduces the order of magnitude of bias.  Call these questions~1 and~2, respectively.  In problem~(ii), it is unclear whether the warp-speed double bootstrap is as effective as the conventional double bootstrap, in the sense of offering the above levels of improved accuracy; we shall refer to this as question~3.  In the present paper we show that the
answers to questions~1 and~2 are positive, but that the answer to question~3 is negative.  In particular, the warp-speed bootstrap does not reduce the order of magnitude of coverage error of a confidence interval.

There is an extensive literature on conventional double-bootstrap methods, particularly in the context of improving the coverage accuracy of single-bootstrap methods.  The first mention of the double bootstrap in this setting apparently was by \cite{Hall_1986}, followed quickly by contributions of \cite{Beran_1987,Beran_1988}.  See also \cite{HallMartin_1988}.  The approach suggested by Hall (1992, Chap.~3) allows general multiple bootstrap methods to be developed together, so that different settings do not require separate treatment.  However, details of properties of the technique seem to be very problem-specific. \cite{Efron_1983} was the first to use the double bootstrap in any setting; in that paper his work was in the context of estimating the error rate of classifiers.  Research on optimising the trade-off between the numbers of simulations in the first and second stages of the conventional double bootstrap, in the context of distribution estimation and constructing confidence
intervals, includes that of \cite{BoothHall_1994}, \cite{BoothPresnell_1998} and \cite{LeeYoung_1999}.

It has become conventional to assess performance of the bootstrap in terms of Edgeworth expansions, not least because that approach enables theoretical properties to be developed in the very broad context addressed by \cite{BhattacharyaGhosh_1978}.  The resulting approximations are valid, in absolute rather than relative terms, uniformly in the tails.  An alternative approach, based on large deviation probabilities, is valid in relative terms; see e.g.~\cite{Hall_1990}.  However, it requires either more stringent assumptions or specialised methods that, at least at present, are not available in the context of the models used by \cite{BhattacharyaGhosh_1978}.  In the setting of absolute rather than relative accuracy, arbitrarily far out into the tails, the results in this paper take the result of consistency, demonstrated by \cite{Giacominietal_2013}, much further.

\section{Model and methodology for bias correction}

\subsection{Model}

Let $\th=f(\mu)$ be a parameter expressible as a known function, $f$, of a $p$-variate mean, $\mu$, and let $\bX$ denote an unbiased estimator of~$\mu=(\mu_1,\ldots,\mu_p)^\T$.  Our estimator of $\th$ is the same function of a sample mean,~$\bX$:
\begin{equation}
\hth=f(\bX)\,.\label{eq:2.1}
\end{equation}
The smooth function $f$ maps a point $x$ in $p$-variate Euclidean space to a point on the real line.  We do not insist that $\bX$ be a mean of $n$, say, independent and identically distributed random $p$-vectors, since it might be the case that $\bX=(\bX_1,\ldots,\bX_p)^\T$, with
$$
\bX_j={1\over n_j}\,\sum_{i=1}^{n_j}\,X\ji\,,
$$
where $X\ji$, for $1\leq i\leq n_j$, are independent for each $i$, $E(X\ji)=\mu_j$ for each $j$, and the $n_j$s are not all equal. Nevertheless, in mathematical terms we shall assume that the $n_j$s are all functions of an integer parameter $n$, and that each $n_j\asymp n$; that is, each ratio $n_j/n$ is bounded away from zero and infinity as $n\rai$.

These issues are related to dependence relationships among the random variables $X\ji$, which should be reflected in resampling methodology.  In our theoretical work we shall suppose that:
\begin{equation} \label{eq:2.2}
\mbox{ \begin{minipage}[c]{.88\linewidth} either (i)~each $n_j=n$ and the vectors $(X_{1i},\ldots,X_{pi})^\T$, for $i\geq1$, are independent and identically distributed; or (ii)~the $X\ji$s are totally independent, for $1\leq i\leq n_j$ and $1\leq j\leq p$, and in this case, for each $j\in\{1,\ldots,p\}$ the variables $X_{j1},X_{j2},\ldots$ are identically distributed, and $n_j\asymp n$.
\end{minipage}}
\end{equation}

Each of (i) and (ii) above can be generalized, for example to hybrid cases where, for positive integers $p_1,\ldots,p_r$ that satisfy $\sum_{j=1}^r p_j=p$, and defining $q_j=\sum_{k=1}^jp_k$, the vectors $V\ji=(X_{q_j+1,i},\ldots,X_{q_{j+1}i})^\T$, for $0\leq j\leq r-1$ and $i\geq1$, are completely independent, and for each $j$ the vectors $V\ji$, for $i\geq1$, are identically distributed. Bootstrap methods that reflect these properties can be constructed readily, and theory providing authoritative support in this setting can be developed, but for the sake of brevity, in our theoretical work we shall restrict attention to cases where (\ref{eq:2.2}) holds.

\subsection{Bias correction}
Bias-corrected estimators of $\th$, based on the conventional bootstrap and the double bootstrap, respectively, are given by
\begin{equation}
\hth\bc=2\,\hth-E(\hth\as\mid\cX)\,,\quad \hth\bcc=3\,\hth-3\,E(\hth\as\mid\cX)+E(\hth\asas\mid\cX)\,.\label{eq:2.3}
\end{equation}
Here $\cX=\{X\ji:1\leq i\leq n_j,1\leq j\leq p\}$ denotes the original dataset, $\hth\as$ is the version of $\hth$ computed from a resample $\cX\as$ drawn randomly, with replacement, from $\cX$, in a manner that reflects appropriately the dependence structure, and $\hth\asas$ is the version of $\hth$ computed from $\cX\asas$, which in turn is drawn randomly with replacement from $\cX\as$, again reflecting dependence.

Monte Carlo approximations to the quantities $\hth\bc$ and $\hth\bcc$ in (\ref{eq:2.3}) are given respectively by
\begin{equation}
\tth\bc=2\,\hth-{1\over B}\,\sum_{b=1}^B\,\hth_b\as\,,\quad \tth\bcc=3\,\hth-{3\over B}\,\sum_{b=1}^B\,\hth_b\as+{1\over BC}\,\sum_{b=1}^B\,\sum_{c=1}^C\,\hth_{bc}\asas\,, \label{eq:2.4}
\end{equation}
where $\hth_b\as$ denotes the $b$th out of $B$ independent and identically distributed, conditional on $\cX$, versions of $\hth\as$, computed from respective resamples $\cX_b\as$ drawn by sampling randomly, with replacement, from the data in $\cX$, and $\hth_{bc}\asas$ is the $c$th out of $C$ independent and identically distributed, conditional on $\cX$ and $\cX\as$, versions of $\hth\asas$, and is computed from a resample $\cX_{bc}\asas$ drawn by sampling randomly, with replacement, from~$\cX_b\as$.

\subsection{Bootstrap algorithms}

Reflecting the model at (\ref{eq:2.1}), we can express $\hth_b\as$ and $\hth_{bc}\asas$ in (\ref{eq:2.4}) as $\hth_b\as=f(\bX_b\as)$ and $\hth_{bc}\asas=f(\bX_{bc}\asas)$, where $\bX_b\as=(\bX_{b1}\as,\ldots,\ab\bX_{bp}\as)^\T$, $\bX_{bc}\asas=(\bX_{bc1}\asas,\ldots,\bX_{bcp}\asas)^\T$, $\bX_{bj}\as$ denotes the mean of data in the resample $\cX_{bj}\as=\{X_{bj1}\as,\ldots,\ab X_{bjn_j}\as\}$, $\bX_{bcj}\asas$ is the mean of data in the re-resample $\cX_{bcj}\asas=\{X_{bcj1}\asas,\ab\ldots,\ab X_{bcjn_j}\asas\}$ drawn by sampling with replacement from $\cX_{bj}\as$, the resampling operations at the first bootstrap level are undertaken by resampling the vectors $X_i=(X_{1i},\ldots,X_{pi})^\T$ randomly, with replacement, if (\ref{eq:2.2})(i) holds, or by resampling the $X\ji$s randomly and completely independently, conditional on $\cX$ and with replacement, if (\ref{eq:2.2})(ii) obtains, and resampling at the second bootstrap level is undertaken analogously.

\subsection{Main conclusions in section~5}

In Theorem~1 in section~5.1 we shall show that if $C\rai$, no matter how slowly, as $n$ and $B$ diverge, then the asymptotic distribution of the Monte Carlo simulation error incurred when constructing $\tth\bcc$ at (\ref{eq:2.4}) is the same as it would be if $C=\infty$.  In particular, not only is the error of order $(nB)\mhf$, the large-sample limiting distribution of the relevant asymptotically normal random variable, which has standard deviation proportional to  $(nB)\mhf$, and which describes in relative detail the accuracy of Monte Carlo bootstrap simulation, is identical to the limiting distribution that would arise if $C=\infty$.

Moreover, if $C$ is held fixed then the order of magnitude, $(nB)\mhf$, remains unchanged, but the standard deviation of the large-sample limiting distribution referred to above changes by a constant factor.  This result is critical.  It demonstrates the relatively small gains that are to be achieved by taking $C$ to be large, and argues in favour of taking $C=1$, for example.  This is the analogue, for bias correction, of the warp-speed bootstrap for distribution estimation when constructing confidence intervals.

Therefore the order of magnitude of Monte Carlo simulation error in $\tth\bcc$ is unchanged even if $C$ is held fixed.  Incidentally, the order of magnitude, $(nB)\mhf$, should be compared with that of the uncorrected bias that remains after applying the bias correction that leads to~$\tth\bcc$; it is~$n\mth$.  Therefore, unless $B$ is of order $n^5$ or larger, for the regular bootstrap, the orders of magnitude involving $B$, discussed above, dominate the error in the bias correction.

\section{ Model and methodology for constructing confidence intervals }

\subsection{Model}

As in section~2.1 we shall assume that the parameter $\th$ can be represented as $f(\mu)$, where the function $f:\IR^p\to\IR$ is known, and $\mu=E(X)$ is an unknown $p$-vector of parameters, estimated by $\bX=n\mo\,\sum_{i=1}^n X_i$ where $\cX=\{X_1,\ldots,X_n\}$ is a random sample of data vectors.  Here and below we use model (\ref{eq:2.2})(i) for the data, but only minor modifications are needed if (\ref{eq:2.2})(ii) is employed instead.

In such cases, provided that $f$ is sufficiently smooth and $\hth$ is given by (\ref{eq:2.1}), the asymptotic variance, $\si_n^2$, of $\hth$ is estimated root-$n$ consistently by $n\mo\,\hsi^2$, where
$$
\hsi^2=\sum_{j_1=1}^p\,\sum_{j_2=1}^p\,f_{j_1j_2}(\bX)\, \oon\,\sumion(X_{j_1i}-\bX_{j_1})(X_{j_2i}-\bX_{j_2})\,.
$$
Here, given a $p$-vector $x=(x_1,\ldots,x_p)^\T$, and integers $j_1,\ldots,j_r$ between 1 and $p$; and assuming that $f$ has $r$ well-defined derivatives with respect to each variable; we put
$$
f_{j_1\ldots j_r}(x)=(\part/\part x_{j_1})\ldots(\part/\part x_{j_r})\,f(x)\,.
$$
The above definitions of $\hth$ and $\hsi$ are used in (\ref{eq:3.1}) below.

\subsection{Bootstrap algorithms}

Let $R$, referred to as the ``root'' by \cite{Giacominietal_2013}, be given by either of the formulae
\begin{equation}
R=n\half\,(\hth-\th)\,,\quad R=n\half\,(\hth-\th)/\hsi\,.\label{eq:3.1}
\end{equation}
Here $\hth$ and $\hsi$ are estimators of parameters $\th$ and $\si$ computed from the random sample $\cX$, and $\si^2$ denotes the asymptotic variance of~$n\half\,\hth$.  The warp-speed bootstrap of \cite{Giacominietal_2013}, closely related to suggestions by \cite{White_2000} and \cite{DavidsonMackinnon_2002,DavidsonMackinnon_2007}, can be defined as follows.

As in section~2, let $\cX_b\as$, for $1\leq b\leq B$, be drawn randomly, with replacement, from $\cX$, and be independent conditional on~$\cX$.  Draw $\cX_b\asas$, denoting a single double-bootstrap resample, by sampling randomly, with replacement, from $\cX_b\as$ for $b=1,\ldots,B$, in such a manner that these re-resamples are independent, conditional on $\cX$ and $\cX_1\as,\ldots,\cX_B\as$.  In the context of section~2, $\cX_b\asas$ would be one of the resamples $\cX_{b1}\asas,\ldots,\cX_{bC}\asas$ which were drawn by resampling from $\cX_b\as$, but on the present occasion we require only one of these resamples.

Let $\hth_b\as$ and $\hth_b\asas$ denote the versions of $\hth$ computed from $\cX_b\as$ and $\cX_b\asas$, respectively, instead of $\cX$, and write $\hsi_b\as$ and $\hsi_b\asas$ for the corresponding versions of~$\hsi$.  If $R$ is given by one of the formulae at (\ref{eq:3.1}), define
\begin{equation}
R_b\as=n\half\,(\hth_b\as-\hth)\,,\quad R_b\as=n\half\,(\hth_b\as-\hth)/\hsi_b\as\,,\label{eq:3.2}
\end{equation}
\begin{equation}
R_b\asas=n\half\,(\hth_b\asas-\hth_b\as)\,,\quad R_b\asas=n\half\,(\hth_b\asas-\hth_b\as)/\hsi_b\asas\,,\label{eq:3.3}
\end{equation}
in the respective cases, and put
\begin{equation}
\hF_B\as(x)=\ooB\,\sum_{b=1}^B\,I(R_b\as\leq x)\,,\quad \tF_B\as(x)=\ooB\,\sum_{b=1}^B\,I(R_b\asas\leq x)\,.\label{eq:3.4}
\end{equation}
Then $\hF_B\as$ is the conventional single-bootstrap, Monte Carlo approximation to the distribution function $F$ of $R$, and the limit of $\hF_B\as$, as $B\rai$, is the conventional single-bootstrap approximation to~$F$.  The function $\tF_B\as$ is a short-cut, warp-speed, double-bootstrap approximation to~$F$.

Given a nominal coverage level $\a\in(0,1)$ of a confidence interval, define $x=\hx_\a\as$ to be the solution of the equation $\tF_B\as(x)=\a$, and similarly let $\hx_\a$ be the solution of $\hF_B\as(x)=\a$.  If $R$ is given by either of the two expressions in (\ref{eq:3.1}), consider the respective confidence intervals,
\begin{equation}
\cI_{b\a}\as=(\hth_b\as-n\mhf\,\hx_\a\as,\infty)\,,\quad \cI_{b\a}\as=(\hth_b\as-n\mhf\,\hsi_b\as\,\hx_\a\as,\infty)\,,\label{eq:3.5}
\end{equation}
which are bootstrap versions of the respective intervals
\begin{equation}
\cI_\a=(\hth-n\mhf\,\hx_\a,\infty)\,,\quad \cI_\a=(\hth-n\mhf\,\hsi\,\hx_\a,\infty)\,.\label{eq:3.6}
\end{equation}
In either case, our estimator of the probability $p_\a$ that the interval $\cI_\a$ covers $\th$ is given by
\begin{equation}
\hp_{B\a}=\ooB\,\sum_{b=1}^B\,I(\hth\in\cI_{b\a}\as)\,.\label{eq:3.7}
\end{equation}
We take the final interval to be $\cI_{\hbe_{B\a}}$, where $\be=\hbe_{B\a}$ denotes the solution of $\hp_{B\be}=\a$.

Earlier warp-speed bootstrap methodology is a little ambiguous in the percentile-$t$ setting, i.e.~in the context of the second definition in each of (\ref{eq:3.1})--(\ref{eq:3.3}), where the technique is not completely clear from the algorithms of \cite{White_2000}, \cite{DavidsonMackinnon_2001,DavidsonMackinnon_2002} and Giacomini \etal(2013, pp.~570--571).  In particular it is unclear from \cite{Giacominietal_2013} when, or whether, the estimator $\hsi$ should be replaced by its single- or double-bootstrap forms, $\hsi\as$ and $\hsi\asas$, for example in (\ref{eq:3.2})--(\ref{eq:3.5}).  The choices we have made are appropriate, however, and in particular the algorithm would not be second-order accurate, or third-order accurate in the case of the double bootstrap, if we were to use simply $\hsi$ in those instances.

\subsection{ Main conclusions drawn in section~5}

In section~5.2 we shall show that in the percentile-$t$ case, using the case $B=\infty$ as a benchmark, the approach suggested above produces quantile estimators that are identical to those obtained using the standard single-bootstrap method, up to an error of order~$n^{-3/2}$.  In particular, they do not reduce the $O(n\mo)$ coverage error of single-bootstrap methods.  Similar results hold for percentile-method bootstrap procedures.

\section{Numerical properties}

\subsection{Bias correction}

Here we report the results of a simulation study comparing the performances of five different bootstrap methods for bias correction: The single bootstrap, the conventional double bootstrap, and the suggested alternative method involving only $C=1$, 2, 5 or 10 double-bootstrap replications.  The data were of two types, either the exponential distribution, with density ${2}^{-1}\,e^{-x/2}$ on the positive half-line, or the log-normal distribution. These two distributions both have nonzero skewness and nonzero kurtosis, making them challenging for the bootstrap. The parameter of interest also took two forms, both of them nonlinear: either $\theta=f(\mu)=\mu^3$ or $\th=\sin(\mu)$, where $\mu$ was the population mean. In such cases there is a term with order $n^{-2}$ in the bias expansion, which cannot be eliminated by the single bootstrap but can be removed by the double bootstrap. This is reflected in our simulation results, which show that the double bootstrap provides better bias
correction than the single bootstrap method.

Sample size, $n$, was chosen in steps of 20 between 20 and~80; the number of simulations, $B$, in the first bootstrap step was set equal to $n^2$, for each of the bootstrap methods; and the number of simulations, $C$, for the second bootstrap step in the conventional double bootstrap was taken to be the integer part of~$10\,B\half$, which we write as $\lfloor 10\,B^{1/2}\rfloor$.  The choice of $B\half$ here was suggested by \cite{BoothHall_1994} in the context of confidence intervals, and gives an expression for $C$ that is orders of magnitude larger than obtained using relatively small, fixed~$C$.  For example, when $n=20$ the value of $C=\lfloor10\,B\half\rfloor$ is between 20 and 200 times the values $C=1$, 2, 5 or 10 used to simulate the alternative approach to double-bootstrap methods; when $n=80$ the respective factors are 80 to~800.

From equation (\ref{eq:2.4}),
\[
\frac{1}{B}\,\sum_{b=1}^B\,\hat{\theta}_b^*-\hat{\theta} ~~\textrm{and}~~\frac{3}{B}\,\sum_{b=1}^B\,\hat{\theta}_b^* -\frac{1}{BC}\,\sum_{b=1}^B\,\sum_{c=1}^C\, \hat{\theta}_{bc}^{**}-2\,\hat{\theta}\,,
\]
provide the estimates of the true bias of $\hat{\theta}$, i.e., $E(\hat{\theta})-\theta$, via single bootstrap and double bootstrap, respectively. Empirical approximations to bias, computed by averaging over the results of 5,000 Monte Carlo trials in each case, are reported in Tables 1-2 in Supplementary Material, and the ratios of such approximations and true bias are graphed in Figure~1.  The figure shows that, for the values of $B$ used in our analysis, there is little to choose between performance when using $C=1$ and $C=\lfloor 10\,B^{1/2}\rfloor$.

\begin{figure}[tbh]
\begin{center}
\subfigure{\includegraphics[scale=0.45]{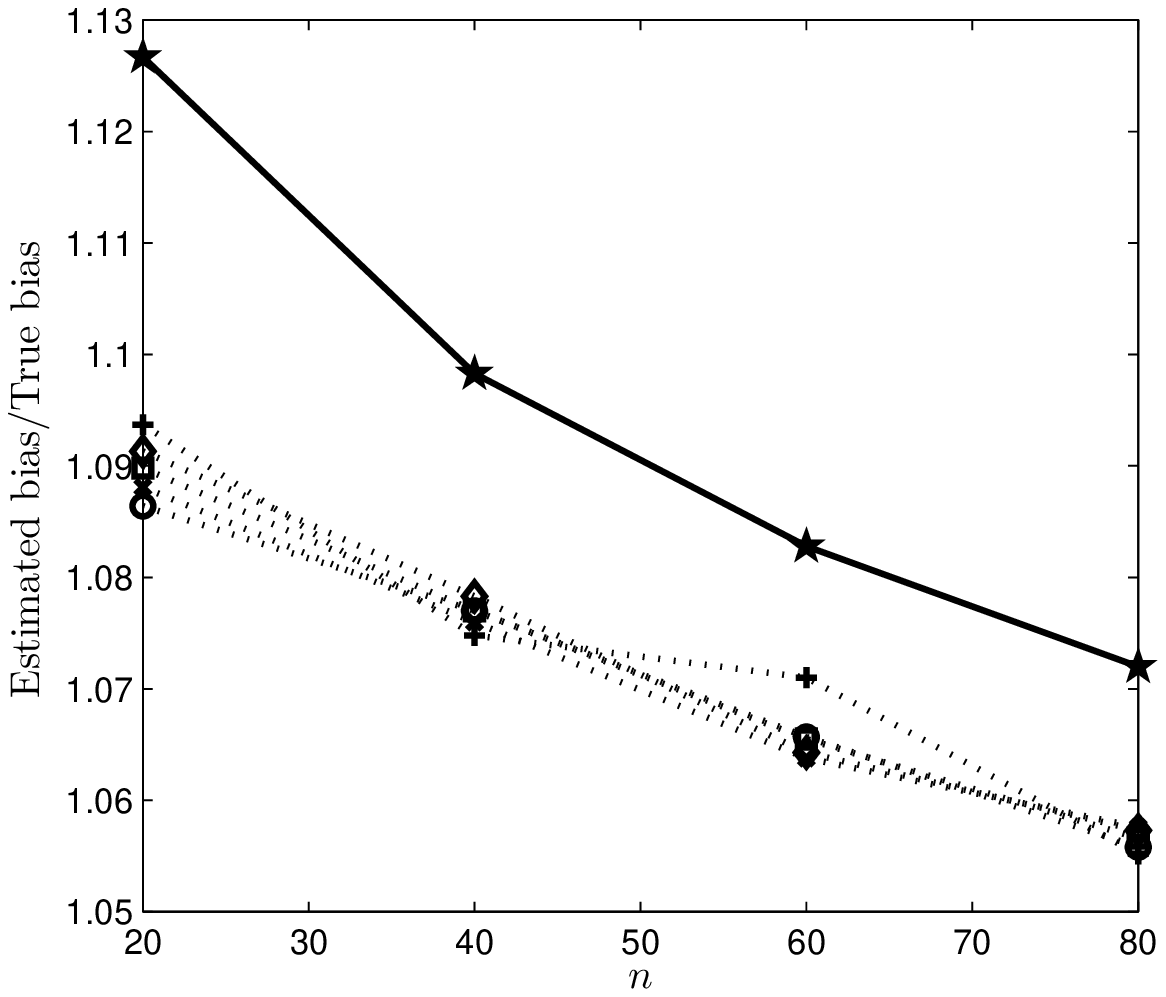}}
\subfigure {\includegraphics[scale=0.45]{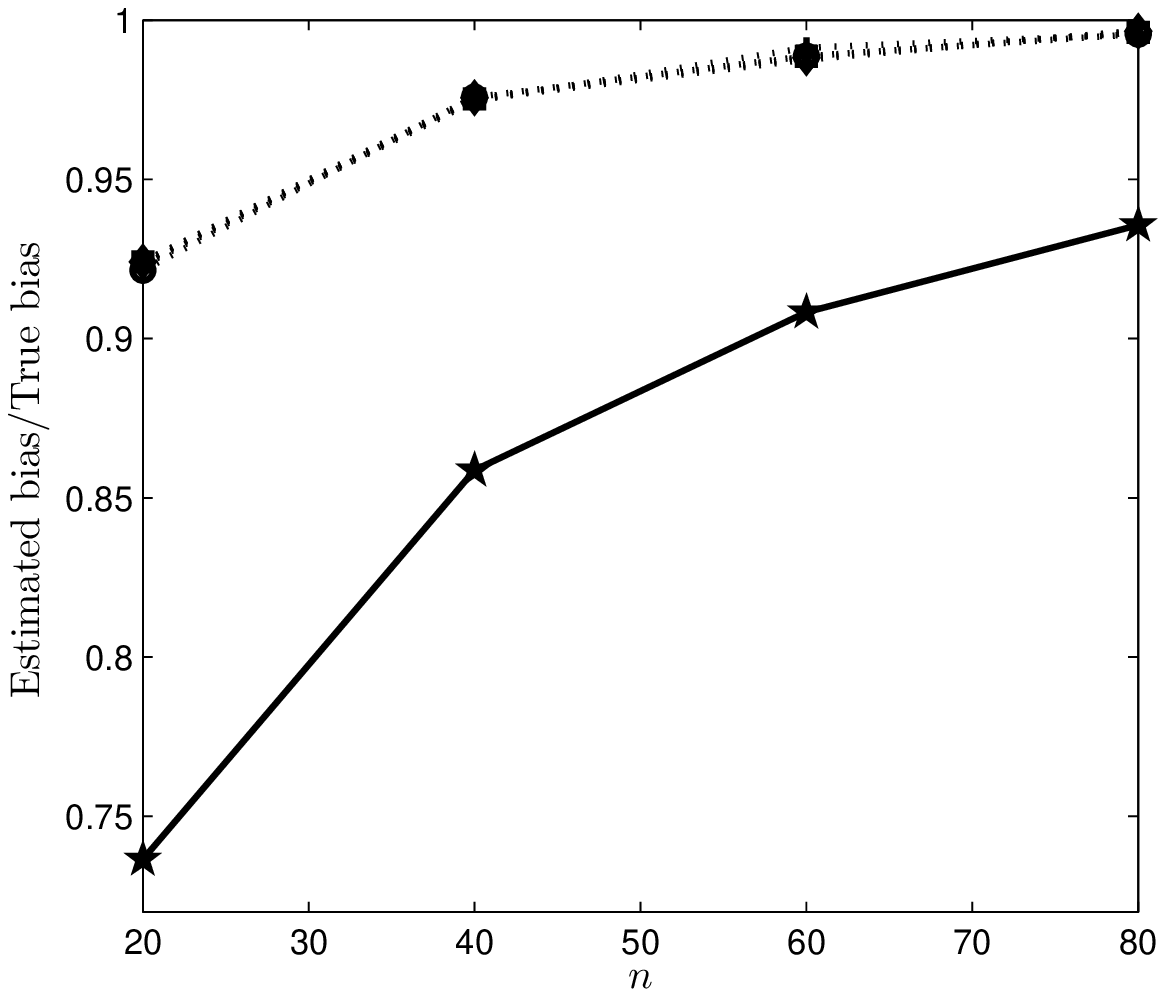}}\\
\subfigure{\includegraphics[scale=0.45]{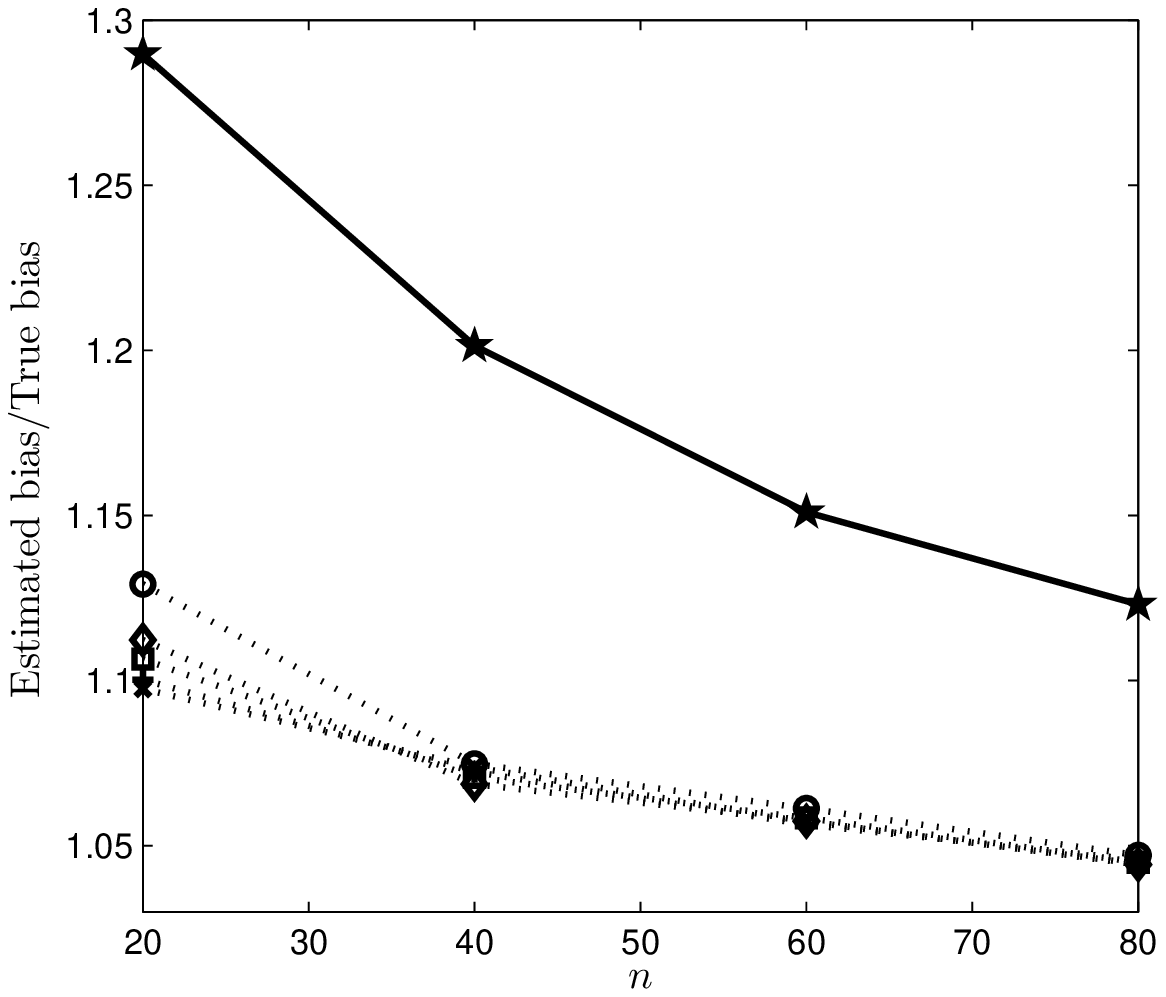}}
\subfigure {\includegraphics[scale=0.45]{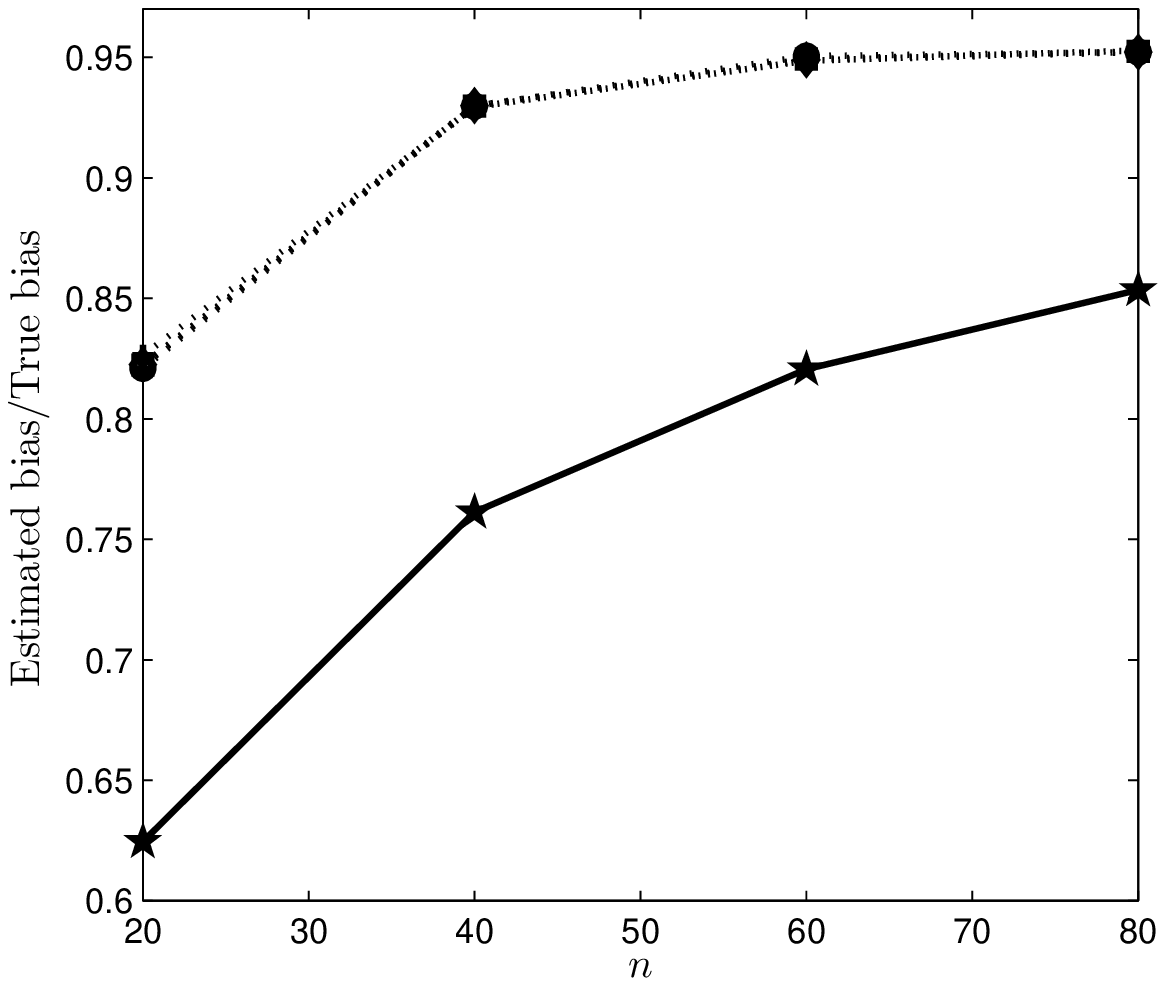}}\\
\end{center}
\caption{ Performance of bootstrap methods for bias correction.  First and second rows show results for the exponential distribution, and the log-normal distribution, respectively; left- and right-hand panels show results for $\th=\mu^3$ and $\th=\sin(\mu)$, respectively. In each panel the graphs represent single bootstrap method ($-\star-$) and conventional double-bootstrap methods with $C=1$ ($\cdots + \cdots$), $C=2$ ($\cdots \circ\cdots$), $C=5$ ($\cdots \times\cdots$), $C=10$ ($\cdots \diamondsuit\cdots$) and $C=\lfloor 10\,B^{1/2}\rfloor$ ($\cdots\square\cdots$), respectively.}
\end{figure}

\subsection{Confidence intervals}

 In this section we illustrate the coverage performance of bootstrap confidence intervals, with nominal coverage $0.9$, for the population means of the two distributions considered in section 4.1, i.e. the exponential and log-normal distributions. Sample size $n$ was taken equal to 20 and 40 in each case; $B$ was increased from $200$ to $700$ in steps of $100$, as indicated on the horizontal axis of each panel; and one-sided and two-sided equal-tailed bootstrap confidence intervals were considered, each using either the percentile or percentile-$t$ bootstrap, implemented via the single bootstrap, the conventional double bootstrap, $C=\lfloor10\,B^{1/2}\rfloor$; and the warp speed bootstrap, i.e.~the double bootstrap with $C=1$. This choice of $C$ was suggested by \cite{LeeYoung_1999}. To provide a perspective different from that in section~4.1, in the present section we graph coverage as a function of $B$ for fixed $n$, rather than as a function of $n$ for fixed $B$ as in section~4.1. Results in the two settings can of course be expressed in same way; the conclusions do not alter.

Results for sample size $n=20$, with each point on each graph based on 5,000 Monte Carlo simulations, are presented in Figure~2. It can be seen that, for each confidence interval type, the conventional double-bootstrap method gives greater coverage accuracy than the single-bootstrap and warp-speed bootstrap. Results for sample size $n=40$ are similar, and are reported in Supplementary Material.

\begin{figure}[tbh]
\begin{center}
\subfigure{\includegraphics[scale=0.45]{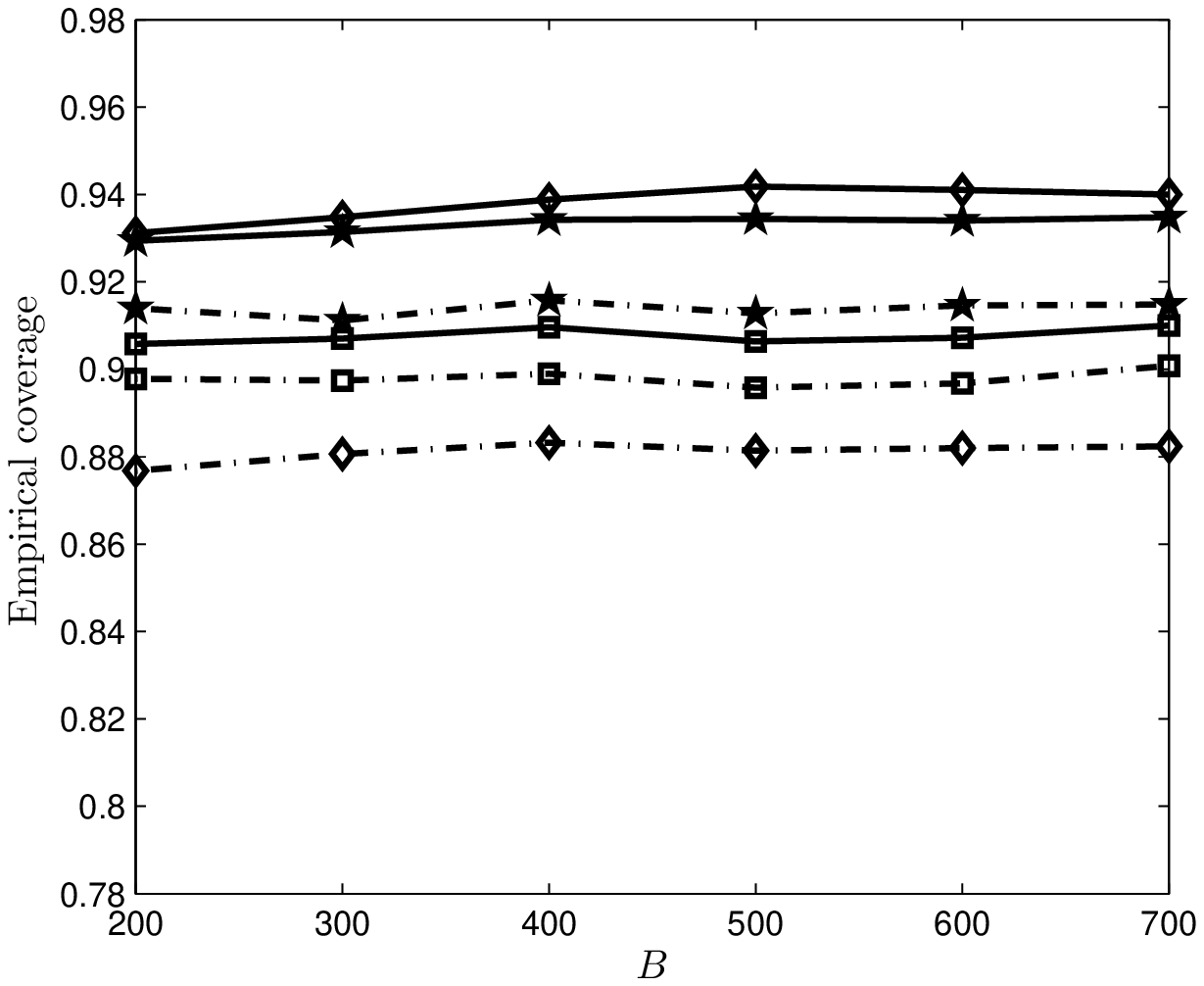}}
\subfigure {\includegraphics[scale=0.45]{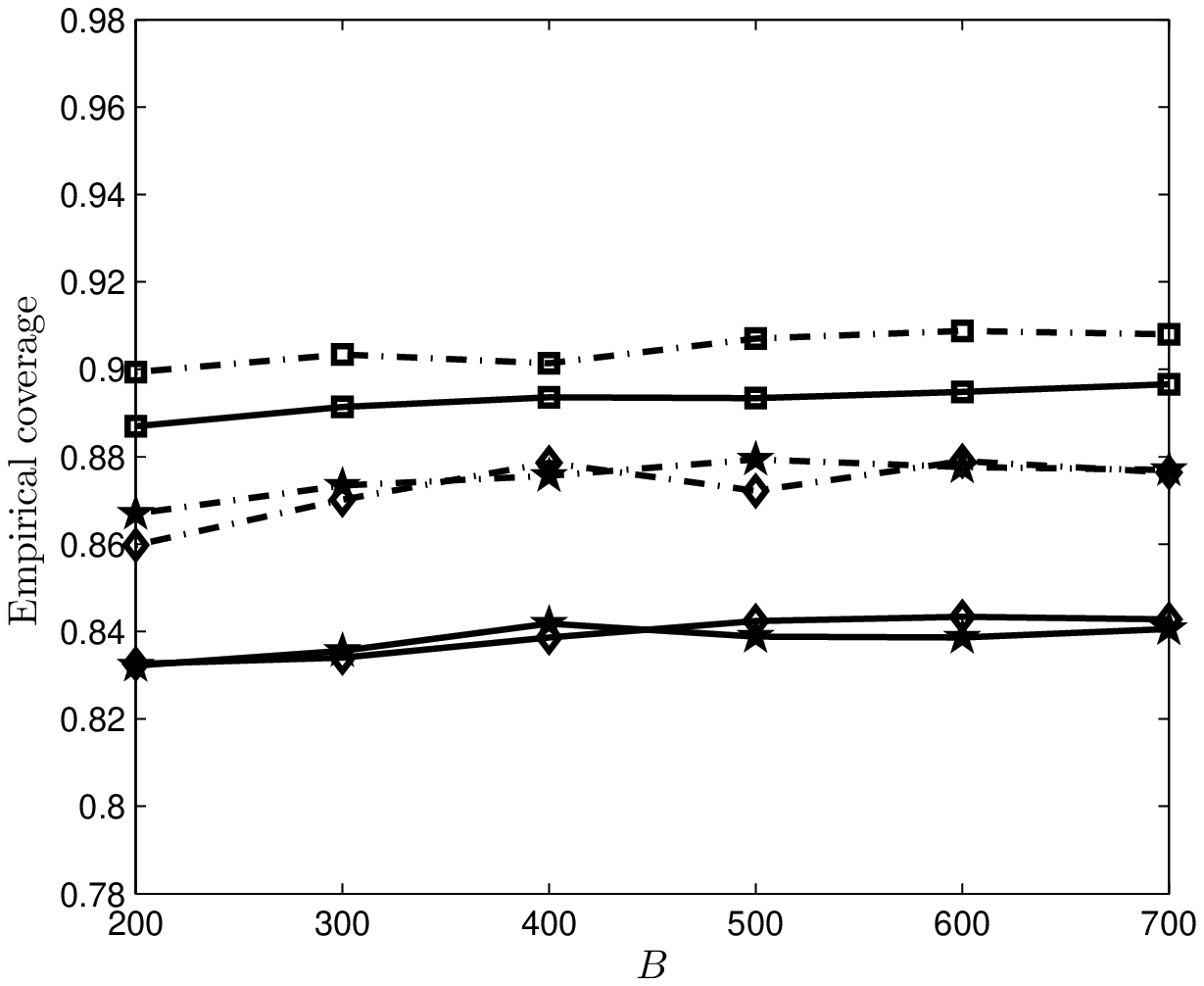}}\\
\subfigure{\includegraphics[scale=0.45]{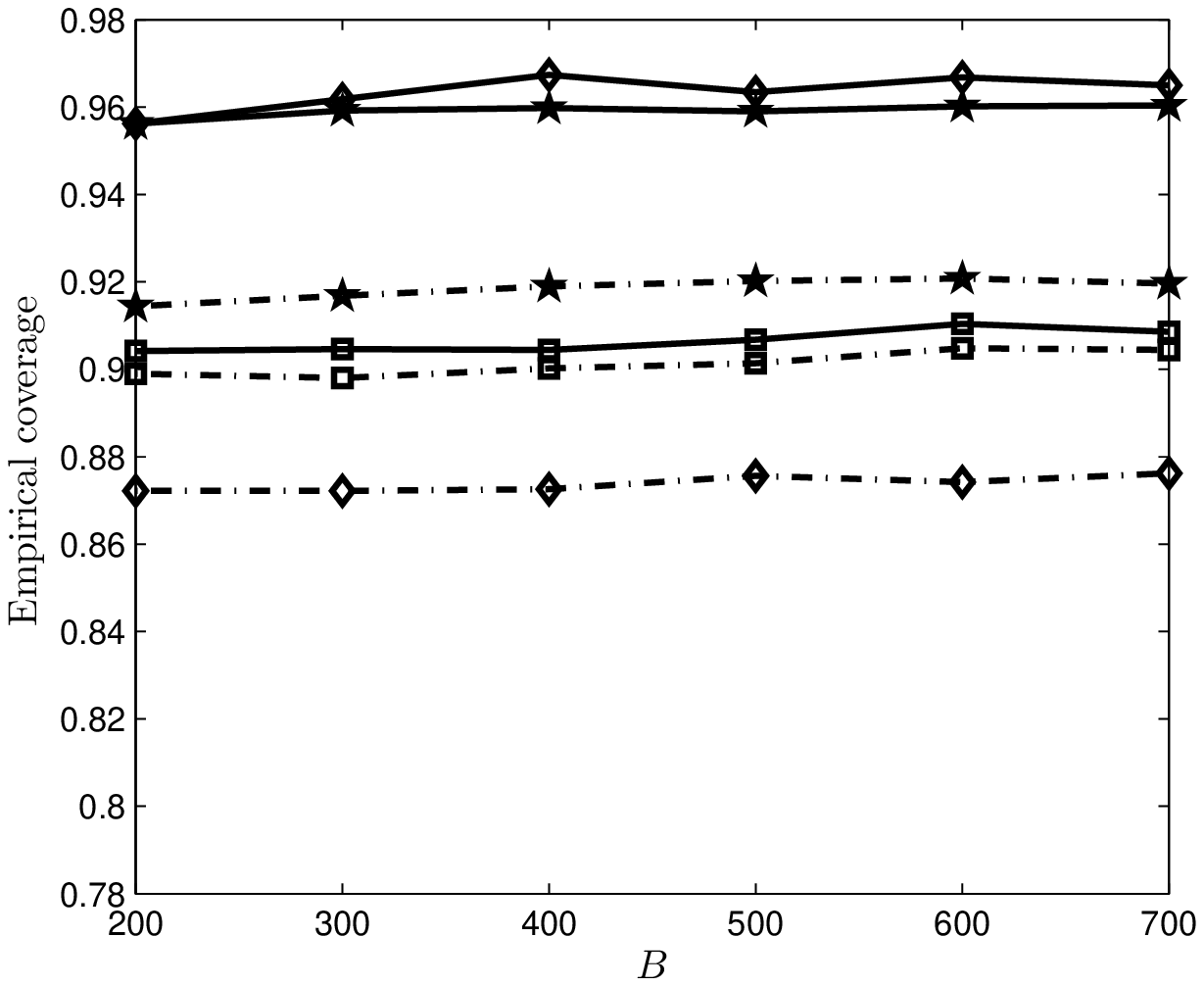}}
\subfigure {\includegraphics[scale=0.45]{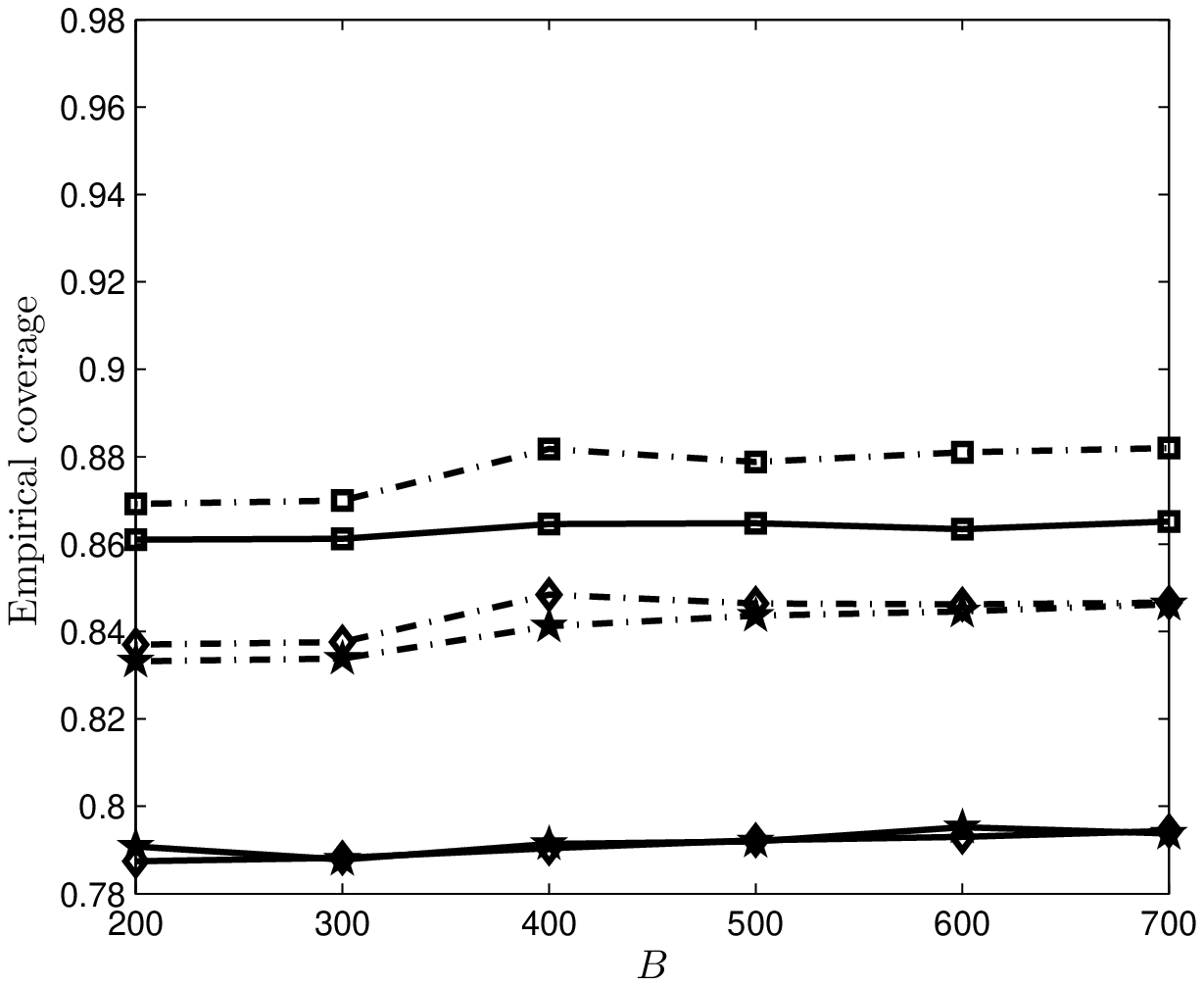}}\\
 \setlength{\abovecaptionskip}{0pt}
\end{center}
\caption{Performance of bootstrap methods for confidence intervals when $n=20$.  First and second rows show results for the exponential distribution, and the log-normal distribution, respectively; left- and right-hand panels show results for one-sided and two-sided equal-tailed confidence intervals, respectively. In each panel the graphs represent single-bootstrap percentile ($-\star-$), single-bootstrap percentile-$t$ ($-\cdot\star\cdot-$), conventional double-bootstrap percentile ($-\square-$), conventional double-bootstrap percentile-$t$ ($-\cdot\square\cdot-$), warp-speed percentile ($-\diamondsuit-$) and warp-speed percentile-$t$ methods ($-\cdot\diamondsuit\cdot-$). }
\end{figure}

\section{Theoretical properties}

\subsection{Bias correction}

Our main regularity condition, in addition to the model assumptions (\ref{eq:2.1}) and (\ref{eq:2.2}), is the following condition:
\begin{equation} \label{eq:5.1}
\mbox{ \begin{minipage}[c]{.88\linewidth} (i) $f(x)$ is differentiable six times with respect to any combination of the $p$ components of $x$; and those derivatives, as well as $f$ itself, are uniformly bounded; and (ii) the data $X_{ji}$ have at least six finite moments, and $E(X_{ji}^6)$ is bounded uniformly in $i$ and $j$.
\end{minipage}}
\end{equation}
Condition (\ref{eq:5.1}) can be generalized, but (for example) if we relax significantly the condition of boundedness of $f$ and its derivatives, in (\ref{eq:5.1})(i), then we need to strengthen the assumption about tails of the distributions of the $X\ji$s, in~(\ref{eq:5.1})(ii).  We shall define
\begin{equation}
\tau^2=E\bigg[\bigg\{\sumjop(X_{j1}-\mu_j)\,f_j(\mu)\bigg\}^{\!2}\bigg]\,.\label{eq:5.2}
\end{equation}

In Theorem~1, below, we decompose the bias-corrected estimators $\tth\bc$, based on the single bootstrap, and $\tth\bcc$, based the double bootstrap, as follows:
\begin{equation}
\tth\bc=U\bc+V\bc\,,\quad\tth\bcc=U\bcc+V\bcc\,,\label{eq:5.3}
\end{equation}
Here $U\bc$ and $U\bcc$ are the ``ideal'' versions of $\tth\bc$ and and $\tth\bcc$, respectively, that we would obtain if we were to do an infinite number of simulations, i.e.~if we were to take $B=C=\infty$; and $V\bc$ and $V\bcc$ denote error terms arising from doing only a finite number of Monte Carlo simulations.  Part~(d) of Theorem~1 shows that the error terms $V\bc$ in the case of the single bootstrap, and $V\bcc$ for the double bootstrap, both equal $O_p\{(nB)\mhf\}$, and that this is the exact order, regardless of the selection of $C$ in the second bootstrap stage. Although the Monte Carlo error terms in the single bootstrap and the double bootstrap share the same convergence rate, equations \eqref{eq:5.4} show that the double bootstrap provides a higher degree of accuracy, in terms of bias correction, than the single bootstrap if we take $B=C=\infty$.  Part~(d) also implies that if $B$ is sufficiently large, or more precisely if $n^5 = O(B)$, then the Monte Carlo error is
of the same order as, or order smaller than, the deterministic remainders in~\eqref{eq:5.4}. These are the main theoretical findings of Theorem~1.

\begin{theorem}Assume that the data are generated according to either of the models at \eqref{eq:2.2}, that \eqref{eq:5.1} holds, and that $B=B(n)\rai$ as $n\rai$.  Then: (a)~Equations \eqref{eq:5.3} hold, where $U\bc$ and $U\bcc$ are functions of $\cX$ alone, and in particular do not involve $\cX\as$ or $\cX\asas$, and satisfy
\begin{equation}
E(U\bc)=\th+O(n\mt)\,,\quad E(U\bcc)=\th+O(n\mth)\,;\label{eq:5.4}
\end{equation}
and $V\bc$ and $V\bcc$ are functions of both $\cX$ and $\cX\as$ (and also of $\cX\asas$, in the case of $V\bcc$), and satisfy $E(V\bc\mi\cX)=E(V\bcc\mi\cX)=0$. (b)~Both $U\bc$ and $U\bcc$ equal $\hth+O_p(n\mo)$, and both satisfy the same central limit theorem as~$\hth$.  (c)~In particular, both $U\bc$ and $U\bcc$ are asymptotically normally distributed with mean $\th$ and a variance, $\si_n^2$ say, which has the property that $n\,\si_n^2$ is bounded as $n\rai$.  (d)~Conditional on $\cX$, $V\bc$ and $V\bcc$ are asymptotically normally distributed with zero means and variances of size $(nB)\mo$, and if $C=C(n)\rai$ as $n\rai$ then the ratio of the variances converges to~1 as $n$ diverges.  In the case of \eqref{eq:2.2}{\rm (i)} the asymptotic variances of $V\bc$ and $V\bcc$, both conditional on $\cX$ and unconditionally, are $(Bn)\mo\,\tau^2$ and $(4+C\mo)\,(Bn)\mo\,\tau^2$, respectively.
\end{theorem}

In connection with part (d) it can be shown that, if $C$ diverges (no matter how slowly) as $n$ increases, the asymptotic distribution of the error is the same as it would be if $C=\infty$. If $\si_n$ is as in part~(c) then, under the model (\ref{eq:2.2})(i), there exists a positive constant $c$ such that $n\,\si_n^2=c+o(1)$ as $n\rai$.  However, this is not necessarily correct under the model (\ref{eq:2.2})(ii), since in that setting we do not require the ratios $n_j/n$ to converge.  In the context of (\ref{eq:2.2})(i), formulae for $U\bc$ and $U\bcc$ are given at (\A9) and (\A10), respectively, in the Supplementary Material.

The orders of magnitude of the remainders in \eqref{eq:5.4} are exact when skewness and kurtosis are nonzero. It follows from part~(b) of Theorem~1 that, in the case $B=C=\infty$, $\tth\bc$ and $\tth\bcc$ satisfy identical central limit theorems, and in particular both have the same asymptotic variances.

\subsection{Distribution estimation and confidence interval construction}

We shall assume that $X$, which represents a generic $p$-vector $X_i=(X_{1i},\ldots,X_{pi})^\T$, where $1\leq i\leq n$ and (\ref{eq:2.2})(i) holds, satisfies the following multivariate version of Cram\'er's continuity condition \citep{Hall_1992}:
\begin{equation}
\limsup_{\|t\|\rai}\,\big|E\{\exp(it^\T X)\}\big|<1\,.\label{eq:5.5}
\end{equation}
On this occasion, $i$ denotes $\sqrt{-1}$.  For brevity we shall treat in detail only the percentile-$t$ case, evidenced by the second formula in each of (\ref{eq:3.1})--(\ref{eq:3.3}), and discuss the percentile method briefly below Theorem~2.

Let $\Phi$ and $\phi$ denote the standard normal distribution and density functions, respectively.  Assume that an unknown scalar parameter $\th$ can be written as $\th=f(\mu)$, where $\mu=E(X)$, and that our estimator of $\th$ is $\hth=f(\bX)$, as at (\ref{eq:2.1}), where $\bX=n\mo\,\sum_{i=1}^n X_i$.  Methods of \cite{BhattacharyaGhosh_1978} can be used to prove that, under conventional assumptions such as those in Theorem~2 below,
\begin{eqnarray}
G(x)&\equiv&
\pr\{n\half\,(\hth-\th)/\hsi\leq x\}\nonumber\\
&=&\Phi(x)+\sum_{j=1}^3\,n^{-j/2}\,Q_j(x)\,\phi(x)+n\mt\,A_n(x)\,,\label{eq:5.6}
\end{eqnarray}
where $Q_j$ is a polynomial of degree $3j-1$, and is an even or odd function according as $j$ is odd or even, respectively; and the remainder $A_n(x)$ satisfies
\begin{equation}
\sup_{n\geq1}\,\sup_{-\infty<x<\infty}\,|A_n(x)|<\infty\,.\label{eq:5.7}
\end{equation}
The coefficients of $Q_j$ are rational polynomials in moments of the distribution of $X$.

For simplicity in this section we take $B=\infty$, which is the ideal case where there is no error generated from Monte Carlo approximation.  Inverting the Edgeworth expansion at \eqref{eq:5.6} we obtain a Cornish-Fisher expansion:
\begin{equation}
x_\a\equiv G\mo(\a)=z_\a+n\mhf\,Q_1\cf(z_\a)+n\mo\,Q_2\cf(z_\a)+n^{-3/2}\,Q_3\cf(z_\a)+O(n\mt)\,, \label{eq:5.A8}
\end{equation}
where $z_\a=\Phi\mo(\a)$, the functions $Q_1\cf$, $Q_2\cf$ and $Q_3\cf$ are Cornish-Fisher polynomials and for example are given by $Q_1\cf=-Q_1$ and $Q_2\cf(x)=Q_1(x)\,Q_1'(x)-\thf\,x\,Q_1(x)^2-Q_2(x)$, and the remainder in \eqref{eq:5.A8} is of the stated order, uniformly in $\a\in[a_1,\a_2]$, whenever $0<\a_1<\a_2<1$.

The conventional percentile-$t$ bootstrap estimator of $G$ is $\hG$, defined~by
\begin{equation*}\label{eq:finite}
\hG(x)=\pr\{n\half\,(\hth\as-\hth)/\hsi\as\leq x\mid\cX\}
\end{equation*}
and satisfying an empirical version of the Edgeworth expansion at \eqref{eq:5.A8}:
\begin{equation}
\hx_\a\equiv \hG\mo(\a)=z_\a+n\mhf\,\hQ_1\cf(z_\a)+n\mo\,\hQ_2\cf(z_\a)+n^{-3/2}\,\hQ_3\cf(z_\a)+O_p(n\mt)\,, \label{eq:5.A9}
\end{equation}
where $\hQ_k\cf$ is derived from empirical Edgeworth polynomials $\hQ_1,\ldots,\hQ_k$ in the standard way, discussed below (\ref{eq:5.A8}); and $\hQ_k$ is derived from the Edgeworth polynomial, $Q_k$, on replacing moments of the distribution of $X$, appearing in coefficients of $Q_k$, by the same respective moments of the distribution of $X\as$, conditional on $\cX$, with $X\as$ drawn by sampling, randomly and with replacement, from~$\cX$.  Note too that the coefficients of $\hQ_k$ depend on moments of $X\as$, conditional on $\cX$, through rational polynomials in those conditional moments.

If we knew the sampling distribution of $X$, and wished to construct an upper one-sided confidence interval for $\th$, we would employ the Studentised confidence interval $(\hth-n\mhf\,\hsi\,x_\a,\infty)$, where $x_\a$ is as at \eqref{eq:5.A8}; if we were to use the percentile-$t$ bootstrap method, it would be $(\hth-n\mhf\,\hsi\,\hx_\a,\infty)$, where $\hx_\a$ is as at \eqref{eq:5.A9}; and if we were to employ the warp-speed bootstrap method, it would be $(\hth-n\mhf\,\hsi\,\hx_{\hbe_\a},\infty)$, as discussed in section~3.2, where $\hbe_\a$ denotes the limit, as $B\rai$, of the quantity $\hbe_{B\a}$ introduced there.  However, we shall show in Theorem~2 that $\hx_{\hbe_\a}=\hx_\a+O_p(n^{-3/2})$, and so the endpoints of standard percentile-$t$ and warp-speed bootstrap confidence intervals differ only to order~$n^{-3/2}$.  This signals that conventional arguments, based on Edgeworth expansions, can be used to prove that the standard percentile-$t$ confidence interval, and its
warp-speed bootstrap variant, have identical coverage error up to and including terms of order $n\mo$, and of course that can be done under the assumptions of Theorem~2.  Since, as is well known, the coverage error of the percentile-$t$ interval is genuinely of order $n\mo$ \citep{Hall_1986}, then it follows that the warp-speed bootstrap does not improve on that accuracy.

\begin{theorem}
Assume that model \eqref{eq:2.2}{\rm (i)} applies; that the function $f$, in the definition $\th=f(\mu)$, has five bounded derivatives; and that \eqref{eq:5.5} holds, $E(\|X\|^K)<\infty$ for sufficiently large $K>0$, and $B=\infty$. Then $\hx_{\hbe_\a}=\hx_\a+O_p(n^{-3/2})$.
\end{theorem}

The appropriate number of moments that should be assumed for general Edgeworth or Cornish Fisher expansions, even in relatively simple, non-bootstrap cases, is awkward to determine.  For example, the argument of \cite{BhattacharyaGhosh_1978} requires at least six moments in the case of the Studentised mean, whereas it is known that three moments are sufficient; see e.g.~\cite{Hall_1987}.  Even if we were to develop, in full detail, a proof of Theorem~2 based on the methods of \cite{BhattacharyaGhosh_1978}, the number of moments we would need to assume would be unduly generous, and instead refer to the number as simply~$K$.  We choose not to provide such a detailed development here.  However, the number of derivatives is relatively easy to address, and the theorem provides detail in that respect.

Let
$$
\widetilde{F}^*(x)=\pr\big\{n^{1/2}\big(\hat{\theta}^{**}-\hat{\theta}^*\big)\big/\hat{\sigma}^{**}\leq x\bigmi\cX\big\}\,,
$$
which is the limit of $\widetilde{F}^*_B(x)$, defined in (\ref{eq:3.4}), as $B\rightarrow\infty$. Then $\hat{x}_{\hat{\beta}_\alpha}$ is the solution of $\widetilde{F}^*(x)=\alpha$. Our focus on the case $B=\infty$ deserves comment.  In the early days of the bootstrap, $B=\infty$ was seen as ``the statistical bootstrap method,'' and the case of finite $B$ was interpreted as a Monte Carlo approximation to the bootstrap.  Indeed, taking $B<\infty$ was viewed more as an issue to be addressed in computational or numerical terms, rather than statistical ones.  Reflecting this, for about eight years from the mid 1980s considerable effort was spent developing efficient computational methods for undertaking bootstrap resampling.  However, by the early 1990s computers had become so fast that this area of research had largely disappeared.  This remains the case today; taking $B$ in the thousands, without using numerical devices to increase simulation efficiency, is now the rule rather than the
exception.  The difference between such large values of $B$, and using the mathematical ideal value $B=\infty$, is particularly small.

\section{Conclusion and discussion}

We have investigated the role played by $C$, the number of resamples used in the second bootstrap stage, in double bootstrap methods for bias correction and confidence intervals.  Specifically, we have shown that the double bootstrap is largely insensitive to choice of $C$ in the context of bias correction.  Indeed, double bootstrap methods with fixed $C$ can produce third-order accuracy, much as do conventional double bootstrap methods with diverging~$C$. This result demonstrates the effectiveness, for bias correction, of using the double bootstrap with a single double-bootstrap simulation.  Although existing work shows that the warp-speed double bootstrap $(C=1)$ can improve accuracy in hypothesis testing, there has not been, until now, any theoretical underpinning of its performance in the context of confidence intervals. However, when only a single bootstrap resample is used in the second-bootstrap stage to construct confidence intervals, the order of magnitude of coverage error
is not improved relative to that for the single bootstrap.

\section*{Supplementary material}
Supplementary Material available for theoretical proofs of Theorems 1 and 2, and additional simulation results in sections 4.1 and 4.2.

\newpage

\renewcommand{\theequation}{A\arabic{equation}}
\setcounter{equation}{0}

\begin{center}
{\bf\Large Supplementary material for ``Double-bootstrap methods use a single double-bootstrap simulation"}\\
\bigskip
Jinyuan Chang\qquad Peter Hall\\
Department of Mathematics and Statistics\\
The University of Melbourne, VIC, 3010, Australia
\end{center}

\section*{A~~~Proof of Theorem~1}

In view of (12), Taylor expansion can be used to derive the following formulae:
\begin{equation}
\hth=\th+\sum_{s=1}^5\,{1\over s!}\,\sum_{j_1=1}^p\ldots\sum_{j_s=1}^p\, (\bX_{j_1}-\mu_{j_1})\cdots(\bX_{j_s}-\mu_{j_s})\,f_{j_1\ldots j_s}(\mu)+O_p(n\mth)\,,
\end{equation}
and
\begin{equation}
E(\hth)=\th+\sum_{s=1}^5\,{1\over s!}\,\sum_{j_1=1}^p\ldots\sum_{j_s=1}^p\, E\{(\bX_{j_1}-\mu_{j_1})\cdots(\bX_{j_s}-\mu_{j_s})\}\,f_{j_1\ldots j_s}(\mu)+O(n\mth)\,,
\end{equation}
where the remainder term $R_n$ that is denoted by $O_p(n\mth)$ in (\A1) satisfies $E(R_n)=O(n\mth)$.

Define
\begin{eqnarray}
\xi_{j_1j_2}&=&\cov(X_{j_11},X_{j_21}),\nonumber\\
\xi_{j_1j_2j_2}&=&E\{(X_{j_11}-\mu_{j_1})\,(X_{j_21}-\mu_{j_2})\,(X_{j_31}-\mu_{j_3})\}\,,\nonumber\\
\xi_{j_1j_2j_3j_4}&=&\xi_{j_1j_2}\,\xi_{j_3j_4} +\xi_{j_1j_3}\,\xi_{j_2j_4}+\xi_{j_1j_4}\,\xi_{j_2j_3}\,.\nonumber
\end{eqnarray}
Then, if (2)(i) holds,
\begin{eqnarray*}
E\{(\bX_{j_1}-\mu_{j_1})(\bX_{j_2}-\mu_{j_2})\}
&=&n\mo\,\xi_{j_1j_2}\,,\\
E\{(\bX_{j_1}-\mu_{j_1})(\bX_{j_2}-\mu_{j_2})(\bX_{j_3}-\mu_{j_3})\}
&=&n\mt\,\xi_{j_1j_2j_3}\,,\\
E\{(\bX_{j_1}-\mu_{j_1})(\bX_{j_2}-\mu_{j_2})(\bX_{j_3}-\mu_{j_3})(\bX_{j_4}-\mu_{j_4})\} &=&n\mt\,\xi_{j_1j_2j_3j_4} +O(n^{-3})\,.
\end{eqnarray*}
Hence, by~(\A2),
\begin{eqnarray}
E(\hth)&=&\th+{1\over2n}\,\sum_{j_1=1}^p\,\sum_{j_2=1}^p\,\xi_{j_1j_2}\,f_{j_1j_2}(\mu)
+{1\over6n^2}\,\sum_{j_1=1}^p\,\sum_{j_2=1}^p\,\sum_{j_3=1}^p\,\xi_{j_1j_2j_3}\,f_{j_1j_2j_3}(\mu)\nonumber\\
&&~~~ +{1\over24\,n^2}\,\sum_{j_1=1}^p\ldots\sum_{j_4=1}^p\,\xi_{j_1j_2j_3j_4}\,f_{j_1j_2j_3j_4}(\mu)
+O(n^{-3})\nonumber\\
&=&\th+\thf\,n\mo\,\ga_2+n\mt\,\big(\osx\,\ga_3+\otf\,\ga_4\big)+O(n^{-3})\,,
\end{eqnarray}
where, for $r=2,3,4$,
$$
\ga_r=\sum_{j_1=1}^p\ldots\sum_{j_r=1}^p\,\xi_{j_1\ldots j_r}\,f_{j_1\ldots j_r}(\mu)\,.
$$

If (2)(ii) holds, instead of (2)(i); and if we define $\si_j^2=\xi_{jj}$, and write $I(\cE)$ for the indicator function of an event $\cE$; then the following relations obtain:
\begin{equation*}
\begin{split}
&~~~~~~~~~~~~~~~~~~~~~~~~~E\{(\bX_{j_1}-\mu_{j_1})(\bX_{j_2}-\mu_{j_2})\}
=n_{j_1}\mo\,I(j_1=j_2)\,\si_{j_1}^2\,,\\
&~~~~~~~~~~E\{(\bX_{j_1}-\mu_{j_1})(\bX_{j_2}-\mu_{j_2})(\bX_{j_3}-\mu_{j_3})\} =n_{j_1}\mt\,I(j_1=j_2=j_3)\,\xi_{j_1j_1j_1}\,,\\
\end{split}
\end{equation*}
and
\begin{equation*}
\begin{split}
&~E\{(\bX_{j_1}-\mu_{j_1})(\bX_{j_2}-\mu_{j_2})(\bX_{j_3}-\mu_{j_3})(\bX_{j_4}-\mu_{j_4})\}\\
=&~(n_{j_1}n_{j_3})\mo\,I(j_1=j_2\neq j_3=j_4)\,\xi_{j_1j_1j_3j_3}+(n_{j_1}n_{j_2})\mo\,I(j_1=j_3\neq
j_2=j_4)\,\xi_{j_1j_2j_1j_2}\\
&~+(n_{j_1}n_{j_2})\mo\,I(j_1=j_4\neq j_2=j_3)\,\xi_{j_1j_2j_2j_1} +O(n^{-3})\,.
\end{split}
\end{equation*}
Therefore we can write (\A2) as
\begin{equation}
E(\hth) =\th+n\mo\,\ga^{(1)}+n\mt\,\ga^{(2)}+O(n^{-3})\,,\label{eq:A.4}
\end{equation}
where the quantities $\ga^{(1)}$ and $\ga^{(2)}$ may depend on $n$ but are bounded as $n\rai$.  Property (\A4) is the analogue, in the context of (2)(ii) rather than (2)(i), of~(\A3).

To explore properties of Monte Carlo approximations to the quantities $E(\hth\as\mi\cX)$ and $E(\hth\asas\mi\cX)$ (compare (3) and (4)), observe first that, analogously to~(\A1),
\begin{eqnarray*}
\hth\as&=&f(\bX\as)=\hth+\sum_{r=1}^5\,{1\over r!}\,\sum_{j_1=1}^p\ldots\sum_{j_r=1}^p\, (\bX_{j_1}\as-\bX_{j_1})\cdots(\bX_{j_r}\as-\bX_{j_r})\,f_{j_1\ldots
j_r}(\bX) +O_p(n\mth)\,,\\
\hth\asas&=&f(\bX\asas)=\hth\as+\sum_{r=1}^5\,{1\over r!}\,\sum_{j_1=1}^p\ldots\sum_{j_r=1}^p\, (\bX_{j_1}\asas-\bX_{j_1}\as)\cdots(\bX_{j_r}\asas-\bX_{j_r}\as)\,f_{j_1\ldots j_r}(\bX\as) +O_p(n\mth)\,.
\end{eqnarray*}
Averaging these formulae over bootstrap replicates we obtain the following expansions:
\begin{eqnarray}
S\bc &\equiv& {1\over B}\,\sum_{b=1}^B\,\hth_b\as=\hth+\sum_{r=1}^5\,{1\over r!}\,\sum_{j_1=1}^p\ldots\sum_{j_r=1}^p\, f_{j_1\ldots
j_r}(\bX)\nonumber \\
 &&~~~~~~~~~~~~~~~~~~~~~~~~~~~~~~~~\times {1\over
B}\,\sum_{b=1}^B\,(\bX_{bj_1}\as-\bX_{j_1})\cdots(\bX_{j_r}\as-\bX_{j_r})
+O_p(n\mth)\,,\\
 S\bcc&\equiv&{1\over
BC}\,\sum_{b=1}^B\,\sum_{c=1}^C\,\hth_{bc}\asas ={1\over B}\,\sum_{b=1}^B\,\hth_b\as+\sum_{r=1}^5\,{1\over r!}\,\sum_{j_1=1}^p\ldots\sum_{j_r=1}^p\, {1\over
B}\,\sum_{b=1}^B\,f_{j_1\ldots j_r}(\bX_b\as)\nonumber\\
&&~~~~~~~~~~~~~~~~~~~~~~~~~~~~~~~\times {1\over C}\,\sum_{c=1}^C\,(\bX_{bcj_1}\asas-\bX_{bj_1}\as)\cdots(\bX_{bcj_r}\asas-\bX_{bj_r}\as)+O_p(n\mth)\,.
\end{eqnarray}
In view of (12), the remainder terms $R_n$, say, that are denoted by $O_p(n\mth)$ in (\A5) and (\A6) satisfy $E(R_n)=O(n\mth)$.

Define
\begin{eqnarray*}
\hxi_{j_1j_2}&=&\oon\,\sumion(X_{j_1i}-\bX_{j_1})(X_{j_2i}-\bX_{j_2})\,,\\
\hxi_{j_1j_2j_2}&=&\oon\,\sumion(X_{j_1i}-\bX_{j_1})(X_{j_2i}-\bX_{j_2})(X_{j_3i}-\bX_{j_3})\,,\\
\hxi_{j_1j_2j_3j_4}&=&\hxi_{j_1j_2}\,\hxi_{j_3j_4}
+\hxi_{j_1j_3}\,\hxi_{j_2j_4}+\hxi_{j_1j_4}\,\hxi_{j_2j_3}\,,\\
\heta_r&=&\sum_{j_1=1}^p\ldots\sum_{j_r=1}^p\,\hxi_{j_1\ldots j_r}\,f_{j_1\ldots j_r}(\bX)\,,
\end{eqnarray*}
the latter for $r=2,3,4$.  In the discussion below we shall assume, for the sake of definiteness, that the data are generated by the model (2)(i).  The case of model (2)(ii) is similar.

Suppose first that we use the regular bootstrap, both for resampling $\cX_b\as$ from $\cX$ and for resampling $\cX\asas_{bc}$ from~$\cX_b\as$.  Then the conditional expected values of the non-remainder terms on the right-hand sides of (\A5) and (\A6) satisfy the following identities, respectively:
\begin{eqnarray}
&&E\bigg\{\hth+\sum_{r=1}^5\,{1\over r!}\,\sum_{j_1=1}^p\ldots\sum_{j_r=1}^p\, f_{j_1\ldots j_r}\big(\bX\big)\, {1\over
B}\,\sum_{b=1}^B\,(\bX_{bj_1}\as-\bX_{j_1})\cdots(\bX_{j_r}\as-\bX_{j_r})\biggmi\cX\bigg\}\nonumber\\
&&~~~~~~~~~~~~~~~~~~=\hth+\thf\,n\mo\,\heta_2
+n\mt\,\big(\osx\,\heta_3+\otf\,\heta_4\big)+O_p(n\mth)\,, \\
&&E\bigg\{{1\over B}\,\sum_{b=1}^B\,\hth_b\as+\sum_{r=1}^5\,{1\over r!}\,\sum_{j_1=1}^p\ldots\sum_{j_r=1}^p\, {1\over
B}\,\sum_{b=1}^B\,f_{j_1\ldots j_r}(\bX_b\as)\nonumber\\
&&~~~~~~~~~~~~~~~~~~~~~~~~~~~~~~~~~~~~~~~~~\times {1\over C}\,\sum_{c=1}^C\,(\bX_{bcj_1}\asas-\bX_{bj_1}\as)\cdots(\bX_{bcj_r}\asas-\bX_{bj_r}\as)\biggmi\cX\bigg\}\cr =&&\hth+\thf\,n\mo\,(2-n\mo)\,\heta_2 +\thf\,n\mt\,\big(\heta_3+\thf\,\heta_4\big)+2\,n\mt\,\big(\osx\,\heta_3+\otf\,\heta_4\big)+O_p(n\mth)\,,
\end{eqnarray}
where, as before, the expected values of the $O_p(n\mth)$ remainder terms equal~$O(n\mth)$.

Recall the definitions of $\tth\bc$ and $\tth\bcc$ at (4), and define
\begin{eqnarray*}
U\bc&\equiv& E(\tth\bc\mid\cX)
=2\,\hth-E(S\bc\mid\cX)\,,\\
 U\bcc&\equiv&
E(\tth\bcc\mid\cX) =3\,\{\hth-E(S\bc\mid\cX)\}+E(S\bcc\mid\cX)\,.
\end{eqnarray*}
Then (\A7) and (\A8) imply that $U\bc=U\bc{}'+O_p(n\mt)$ and $U\bcc=U\bcc{}'+O_p(n\mth)$, where the expected values of the $O_p(n^{-k})$ remainder terms equal $O(n^{-k})$, and
\begin{eqnarray}
 U\bc{}'&=&\hth-\big\{\thf\,n\mo\,\heta_2
+n\mt\,\big(\osx\,\heta_3+\otf\,\heta_4\big)\big\}\,,\\
U\bcc{}'&=&\hth-\thf\,n\mo\,(1+n\mo)\,\heta_2 +\thf\,n\mt\,\big(\heta_3+\thf\,\heta_4\big) -n\mt\,\big(\osx\,\heta_3+\otf\,\heta_4\big)\,.\quad\;
\end{eqnarray}
Therefore $U\bc$ and $U\bcc$ both equal $\hth+O_p(n\mo)$, as claimed in part~(b) of Theorem~1.

Put $V\bc=\tth\bc-E(\tth\bc\mi\cX)$ and $V\bcc=\tth\bcc-E(\tth\bcc\mi\cX)$.  Employing (\A3) and the properties
\begin{equation}
E(\heta_2)=(1-n\mo)\,\ga_2+n\mo\,\big(\ga_3+\thf\,\ga_4\big)+O(n\mt)\,,\quad E(\heta_r)=\ga_r+O(n\mo)\;\label{eq:A.11}
\end{equation}
for $r=3,4$, we deduce that $E(\tth\bc)=E(U\bc)=\th+O(n\mt)$, and that $V\bc=\tth\bc-U\bc$ is a function of both $\cX$ and $\cX\as$, satisfying $E(V\bc\mi\cX)=0$ (in the context of (2)(i)) and $\var(V\bc\mi\cX)=\{1+o_p(1)\}\,(Bn)\mo\,\tau^2$.  Central limit theorems for $U\bc$ and $V\bc$ follow from Lindeberg's theorem.  In the context of (2)(i), those parts of (15) and (b)--(d), in Theorem~1, that pertain to the single-bootstrap estimator $\tth\bc$, follow from these properties.  (The exactness of the orders of magnitude of remainders in (15) can be proved by deriving concise formulae for those terms, using (\A9)--(\A11).)

The results discussed two paragraphs above also imply that $E(\tth\bcc)=E(U\bcc)\ab=\th+O(n\mth)$, and of course, $V\bcc=\tth\bcc-U\bcc$ is a function of $\cX$, $\cX\as$ and $\cX\asas$ satisfying $E(V\bcc\mi\cX)=0$.  Note too that, in the context of (2)(i),
\begin{eqnarray*}
&&(BC)^2\,\var(S\bcc-S\bc\bigmi\cX) \sim_p\var\bigg\{\sum_{b=1}^B\,\sum_{c=1}^C\,\sumjop
f_j(\bX_b\as)\,(\bX_{bcj}\asas-\bX_{bj}\as)\biggmi\cX\bigg\}\\
&&~~~~~~~~~~~~~~~~~~~=E\Bigg[\bigg\{\sum_{b=1}^B\,\sum_{c=1}^C\,\sumjop
f_j(\bX_b\as)\,(\bX_{bcj}\asas-\bX_{bj}\as)\bigg\}^{\!2}\Biggmi\cX\Bigg]\\
&&~~~~~~~~~~~~~~~~~~~\sim_pE\Bigg[\bigg\{\sum_{b=1}^B\,\sum_{c=1}^C\,\sumjop
f_j(\mu)\,(\bX_{bcj}\asas-\bX_{bj}\as)\bigg\}^{\!2}\Biggmi\cX\Bigg]\\
&&~~~~~~~~~~~~~~~~~~~=E\Bigg(E\Bigg[\bigg\{\sum_{b=1}^B\,\sum_{c=1}^C\,\sumjop
f_j(\mu)\,(\bX_{bcj}\asas-\bX_{bj}\as)\bigg\}^{\!2}\Biggmi\cX,\cX\as\Bigg]\Biggmi\cX\Bigg)\\
&&~~~~~~~~~~~~~~~~~~~=C\,E\Bigg(E\Bigg[\bigg\{\sum_{b=1}^B\,\sumjop
f_j(\mu)\,(\bX_{b1j}\asas-\bX_{bj}\as)\bigg\}^{\!2}\Biggmi\cX,\cX\as\Bigg]\Biggmi\cX\Bigg)\\
&&~~~~~~~~~~~~~~~~~~~=C\,E\Bigg[\bigg\{\sum_{b=1}^B\,\sumjop
f_j(\mu)\,(\bX_{b1j}\asas-\bX_{bj}\as)\bigg\}^{\!2}\Biggmi\cX\Bigg]\\
&&~~~~~~~~~~~~~~~~~~~=BC\,E\Bigg[\bigg\{\sumjop f_j(\mu)\,(\bX_{11j}\asas-\bX_{1j}\as)\bigg\}^{\!2}\Biggmi\cX\Bigg] \sim_p BC\,n\mo\,\tau^2\,,
\end{eqnarray*} and
$\cov(S\bcc-S\bc,S\bc\mi\cX)=o_p(B\mo)$.  Therefore,
\begin{equation*}
\begin{split}
\var(\tth\bcc\mid\cX) &=\var(V\bcc\mid\cX)
=\var(S\bcc-3\,S\bc\mid\cX)\\
&=\var(S\bcc-S\bc\mid\cX)
-4\,\cov(S\bcc-S\bc,S\bc\mid\cX)+4\,\var(S\bc\mid\cX)\\
&=(nB)\mo\,(4+C\mo)\,\tau^2+o_p\{(nB)\mo\}\,.
\end{split}
\end{equation*}
Much as in the case of $\tth\bc$, it can be proved from (\A10) and (\A11) that $E(\tth\bcc)=E(U\bcc)=\th+O(n\mth)$.  If (2)(i) holds then these properties, and Lindeberg's central limit theorem, imply those parts of Theorem~1 that pertain to the double-bootstrap estimator $\tth\bcc$.   Cases where the model (2)(ii) holds are similar.

\section*{B~~~Proof of Theorem~2}

Consider first the solution $\be=\be_\a$, say, of the equation
\begin{equation}
\pr\{n\half\,(\hth\as-\hth)/\hsi\as\leq x_\be\} =\a\,,\label{eq:A.A12}
\end{equation}
where $x=x_\be$ is the solution of
\begin{equation}
\pr\{n\half\,(\hth-\th)/\hsi\leq x\}=\be\,.\label{eq:A.A13}
\end{equation}
Note that
\begin{equation}
\begin{split}
&~\pr\{n\half\,(\hth\as-\hth)/\hsi\as\leq x\mid\cX\}\\
=&~\Phi(x)+n\mhf\,\hQ_1(x)\,\phi(x)+\cdots+n^{-3/2}\,\hQ_3(x)\,\phi(x) +n\mt\,\hA_{n1}(x)\,,
\end{split}\label{eq:A.A14}
\end{equation}
where the remainder $\hA_{n1}(x)$ satisfies
\begin{equation}
\pr\bigg\{\sup_{-\infty<x<\infty}|\hA_{n1}(x)|>n^{K_1}\bigg\}=O(n^{-K_2})\label{eq:A.A15}
\end{equation}
and the constants $K_1$ and $K_2$, both of which are strictly positive, can be chosen as small or as large, respectively, as desired, at the expense of having to assume a higher moment of $\|X\|$ in the theorem.

The left-hand side of \eqref{eq:A.A12} equals the expected value of the left-hand side of \eqref{eq:A.A14}, and hence also of the right-hand side of that formula.  The coefficients of $\hQ_k$ depend on moments of $X\as$, conditional on $\cX$, through rational polynomials in those conditional moments.  The denominators in those rational polynomials can be Taylor expanded, obtaining quantities $\hQ_k^{{\rm exp}}$, say, which have the property that
$$
\sup_{-\infty<x<\infty}\Big[E\big\{\big|\hQ_k^{{\rm exp}}(x)\big|\big\}\,\phi(x)\Big]=O(1)\,,\quad E\big\{\hQ_k^{{\rm exp}}(x)\big\}\,\phi(x)=Q_k(x)\,\phi(x)+O(n\mo)\,,
$$
where the latter identity holds uniformly in $x$; and also,
\begin{equation*}
\begin{split}
n\mhf\,\hQ_1(x)\,\phi(x)&+\cdots+n^{-3/2}\,\hQ_3(x)\,\phi(x)
+n\mt\,\hA_{n1}(x)\\
&=n\mhf\,\hQ_1^{{\rm exp}}(x)\,\phi(x)+\cdots+n^{-3/2}\,\hQ_3^{{\rm exp}}(x)\,\phi(x) +n\mt\,\hA_{n2}(x)\,,
\end{split}
\end{equation*}
and $\hA_{n2}$ satisfies \eqref{eq:A.A15} and additionally, $E\{\hA_{n2}(x)\}=O(n\half)$, uniformly in~$x$.  Hence, taking the expected value of both sides of \eqref{eq:A.A14}, we deduce that
\begin{equation}
\begin{split}
&~\pr\{n\half\,(\hth\as-\hth)/\hsi\as\leq x\}\\
=&~\Phi(x)+n\mhf\,Q_1(x)\,\phi(x)+\cdots+n^{-3/2}\,Q_3(x)\,\phi(x) +O(n^{-3/2})\,,
\end{split}\label{eq:A.A16}
\end{equation}
from which it follows that
$$
\pr\{n\half\,(\hth\as-\hth)/\hsi\as\leq x\} =\pr\{n\half\,(\hth-\th)/\hsi\leq x\} +O(n^{-3/2})\,.
$$

However, the solution $x=x_{\be_\a}$ of
$$
\pr\{n\half\,(\hth-\th)/\hsi\leq x\} +O(n^{-3/2}) =\a
$$
is identical, up to terms of order $n^{-3/2}$, to the solution $x=x_\a$ of equation \eqref{eq:A.A13} when $\be=\a$ there, and in particular,
$$
x_{\be_\a}=x_\a+O(n^{-3/2})\,.
$$
Therefore,
\begin{equation}
x_{\be_\a}=z_\a+n\mhf\,Q_1\cf(z_\a)+n\mo\,Q_2\cf(z_\a)+O(n^{-3/2})\,.\label{eq:A.A17}
\end{equation}

Recall that the distribution function estimator with which we are working is the version of the second formula in (8) when $B=\infty$ and $C=1$:
$$
\tF_\infty(x)=\pr\{n\half\,(\hth\asas-\hth\as)/\hsi\asas\leq x\mid\cX\}\,,
$$
where $\hth\as$, $\hth\asas$ and $\hsi\asas$ are computed from $\cX\as$, $\cX\asas$ and $\cX\asas$, respectively.  Since we are taking $B=\infty$ in our analysis then $\hx_\a$, defined (9) in the case of finite $B$, is now given by the limit as $B\rai$ of that definition, i.e.~the solution in $x$ of $\pr\{n\half\,(\hth\as-\hth)/\hsi\as\leq x\mi\cX\}=\a$.  In this notation, $\hbe_\a$ is defined to be the solution in $\be$ of the equation $\tF_\infty(\hx_\be)=\a$, i.e.~the solution in $\be$ of
\begin{equation}
\pr\{n\half\,(\hth\asas-\hth\as)/\hsi\asas\leq\hx_\be\mid\cX\} =\a\,.\label{eq:A.A18}
\end{equation}
Now, the solution in $\be$ of \eqref{eq:A.A18} is an estimator of the solution $\be=\be_\a$ of
$$
\pr\{n\half\,(\hth\as-\hth)/\hsi\as\leq x_\be\} =\a\,,
$$
where $x=x_\be$ is the solution of \eqref{eq:A.A13}.  That is, a representation of $\hx_{\hbe_\a}$ as a Cornish-Fisher expansion is identical to the analogous representation of $x_{\be_\a}$, except that moments of $X$ are replaced by the corresponding moments of $X\as$ conditional on~$\cX$.  Since the Cornish-Fisher expansion of $x_{\be_\a}$ is given by \eqref{eq:A.A17}, up to and including terms of order $n\mo$, then
$$
\hx_{\hbe_\a}=z_\a+n\mhf\,\hQ_1\cf(z_\a)+n\mo\,\hQ_2\cf(z_\a)+O_p(n^{-3/2})\,.
$$
This is identical to the expansion of $\hx_\a$, the solution of
$$
\pr\{n\half\,(\hth\as-\hth)/\hsi\as\leq x\mid\cX\}=\a\,,
$$
up to and including terms of order $n\mo$, and so $\hx_{\hbe_\a}=\hx_\a+O_p(n^{-3/2})$, as had to be proved.

\section*{C~~~Simulation results}

In this section, we provide the simulation results for sections 4.1 and 4.2.
\subsection*{C.1~~~Bias estimation in section 4.1}

Tables 1 and 2 report the empirical approximations to bias computed by averaging over the results of 5,000 Monte Carlo trials in the settings of exponential distribution and log-normal distribution, respectively.
\begin{table}
\caption{Bias estimation  based on different bootstrap methods for $\mu^3$ and $\sin(\mu)$ with Exp(2) distribution. The values in brackets denote the ratios of the estimated biases and true bias, respectively.}
\bigskip
\centering{  \begin{tabular}{llccccc}
     &  $n$&$20$  & $40$  & $60$  &  $80$ \\[5pt]
    $\mu^3$&true bias $(\times 10^2)$                                                   & 115.1658 & 57.0163 & 38.1427 & 28.6419  \\
           &        single $(\times 10^2)$                                       & 129.7612 & 62.6221 & 41.3012 & 30.7055 \\
           &                                                                     & [1.1267] & [1.0983]& [1.0828]& [1.0720]\\
           &        double with $C=1$  $(\times 10^2)$                           & 125.9539 & 61.2805 & 40.8512 & 30.2225 \\
           &                                                                     & [1.0937] & [1.0748]& [1.0710]& [1.0552]\\
           &        double with $C=2$  $(\times 10^2)$                           & 125.1125 & 61.4080 & 40.6490 & 30.2391\\
           &                                                                     & [1.0864] & [1.0770]& [1.0657]& [1.0558]\\
           &        double with $C=5$   $(\times 10^2)$                          & 125.3128 & 61.3515 & 40.5743 & 30.2928\\
           &                                                                     & [1.0881] & [1.0760]& [1.0638]& [1.0576]\\
           &        double with $C=10$   $(\times 10^2)$                         & 125.6812 & 61.4801 & 40.5936 & 30.2841\\
           &                                                                     & [1.0913] & [1.0783]& [1.0643]& [1.0573]\\
           &        double with $C=\lfloor 10B^{1/2} \rfloor$ $(\times 10^2)$    & 125.5125 & 61.4068 & 40.6418 & 30.2630  \\
           &                                                                     & [1.0898] & [1.0770]& [1.0655]& [1.0566]\\[5pt]
    $\sin(\mu)$ &true bias $(\times 10^2)$                                       & -8.4970 & -4.4585 & -2.9896 & -2.2458  \\
            &      single $(\times 10^2)$                                        & -6.2578 & -3.8283 & -2.7155 & -2.1012\\
            &                                                                    & [0.7365]& [0.8587]& [0.9083]& [0.9356]\\
            &      double with $C=1$  $(\times 10^2)$                            & -7.8440 & -4.3452 & -2.9636 & -2.2358\\
            &                                                                    & [0.9231]& [0.9746]& [0.9913]& [0.9955]\\
            &      double with $C=2$  $(\times 10^2)$                            & -7.8299 & -4.3505 & -2.9557 & -2.2359\\
            &                                                                    & [0.9215]& [0.9758]& [0.9887]& [0.9956]\\
            &      double with $C=5$   $(\times 10^2)$                           & -7.8483 & -4.3475 & -2.9526 & -2.2383\\
            &                                                                    & [0.9237]& [0.9751]& [0.9876]& [0.9967]\\
            &      double with $C=10$   $(\times 10^2)$                          & -7.8521 & -4.3499 & -2.9541 & -2.2380\\
            &                                                                    & [0.9241]& [0.9756]& [0.9881]& [0.9965]\\
            &      double with $C=\lfloor 2B^{1/2} \rfloor$ $(\times 10^2)$      & -7.8520 & -4.3480 & -2.9555 & -2.2371  \\
            &                                                                    & [0.9241]& [0.9752]& [0.9886]& [0.9961]\\
\end{tabular}}
\end{table}

\begin{table}
\caption{Bias estimation  based on different bootstrap methods for $\mu^3$ and $\sin(\mu)$ with $\exp\{N(0,1)\}$ distribution. The values in brackets denote the ratios of the estimated biases and true bias, respectively.}
\bigskip
\centering{  \begin{tabular}{llccccc}
      & $n$&$20$  & $40$  & $60$  &  $80$ \\[5pt]
    $\mu^3$ &true bias $(\times 10^2)$                                                   & 116.4471 & 55.6341 & 36.9453 & 27.9352  \\
            &       single $(\times 10^2)$                                        & 150.1797 & 66.8400 & 42.5223 & 31.3730\\
            &                                                                     & [1.2897] & [1.2014]& [1.1510]& [1.1231]\\
            &       double with $C=1$  $(\times 10^2)$                            & 128.1239 & 59.6595 & 39.0303 & 29.2126\\
            &                                                                     & [1.1003] & [1.0724]& [1.0564]& [1.0457]\\
            &       double with $C=2$  $(\times 10^2)$                            & 131.4972 & 59.7961 & 39.2092 & 29.2521\\
            &                                                                     & [1.1292] & [1.0748]& [1.0613]& [1.0471]\\
            &       double with $C=5$   $(\times 10^2)$                           & 127.7990 & 59.7409 & 39.0654 & 29.1772\\
            &                                                                     & [1.0975] & [1.0738]& [1.0574]& [1.0445]\\
            &       double with $C=10$   $(\times 10^2)$                          & 129.5233 & 59.4563 & 39.0700 & 29.1729\\
            &                                                                    & [1.1123] & [1.0687]& [1.0575]& [1.0443]\\
            &      double with $C=\lfloor 10B^{1/2} \rfloor$ $(\times 10^2)$     & 128.8509 & 59.5656 & 39.1011 & 29.1925  \\
            &                                                                     & [1.1065] & [1.0707]& [1.0584]& [1.0450]\\[5pt]
    $\sin(\mu)$ &true bias $(\times 10^2)$                                        & -9.8256 & -5.6652 & -3.9217 &  -2.9741 \\
     &              single $(\times 10^2)$                                        & -6.1373 & -4.3128 & -3.2181 & -2.5383\\
     &                                                                            & [0.6246]& [0.7613]& [0.8206]& [0.8535]\\
     &              double with $C=1$  $(\times 10^2)$                            & -8.1200 & -5.2653 & -3.7202 & -2.8340\\
     &                                                                            & [0.8264]& [0.9294]& [0.9486]& [0.9529]\\
     &              double with $C=2$  $(\times 10^2)$                            & -8.0672 & -5.2670 & -3.7275 & -2.8318\\
     &                                                                            & [0.8210]& [0.9297]& [0.9505]& [0.9522]\\
     &              double with $C=5$   $(\times 10^2)$                           & -8.0785 & -5.2651 & -3.7201 & -2.8321\\
     &                                                                            & [0.8222]& [0.9294]& [0.9486]& [0.9523]\\
     &              double with $C=10$   $(\times 10^2)$                          & -8.0812 & -5.2684 & -3.7214 & -2.8320\\
     &                                                                            & [0.8225]& [0.9300]& [0.9489]& [0.9522]\\
     &              double with $C=\lfloor 10B^{1/2} \rfloor$ $(\times 10^2)$     & -8.0796 & -5.2667 & -3.7228 & -2.8324  \\
     &                                                                            & [0.8223]& [0.9297]& [0.9493]& [0.9524]\\
\end{tabular}}
\end{table}

\subsection*{C.2~~~Performance of $n=40$ in section 4.2}

Figure 3 shows the empirical coverage of the confidence intervals constructed by different bootstrap methods when sample size $n=40$.
\begin{figure}[tbh]
\begin{center}
\subfigure{\includegraphics[scale=0.45]{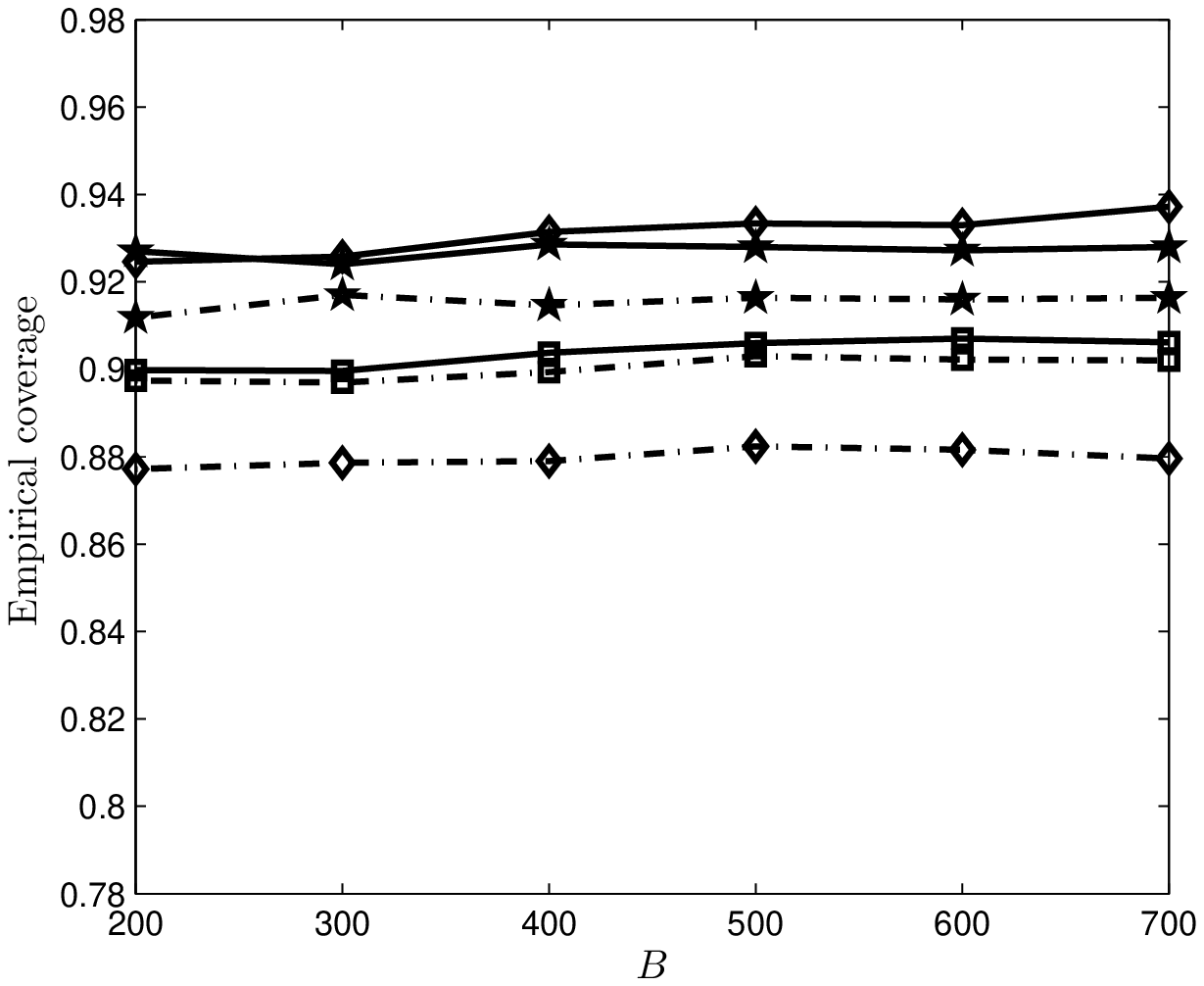}}
\subfigure {\includegraphics[scale=0.45]{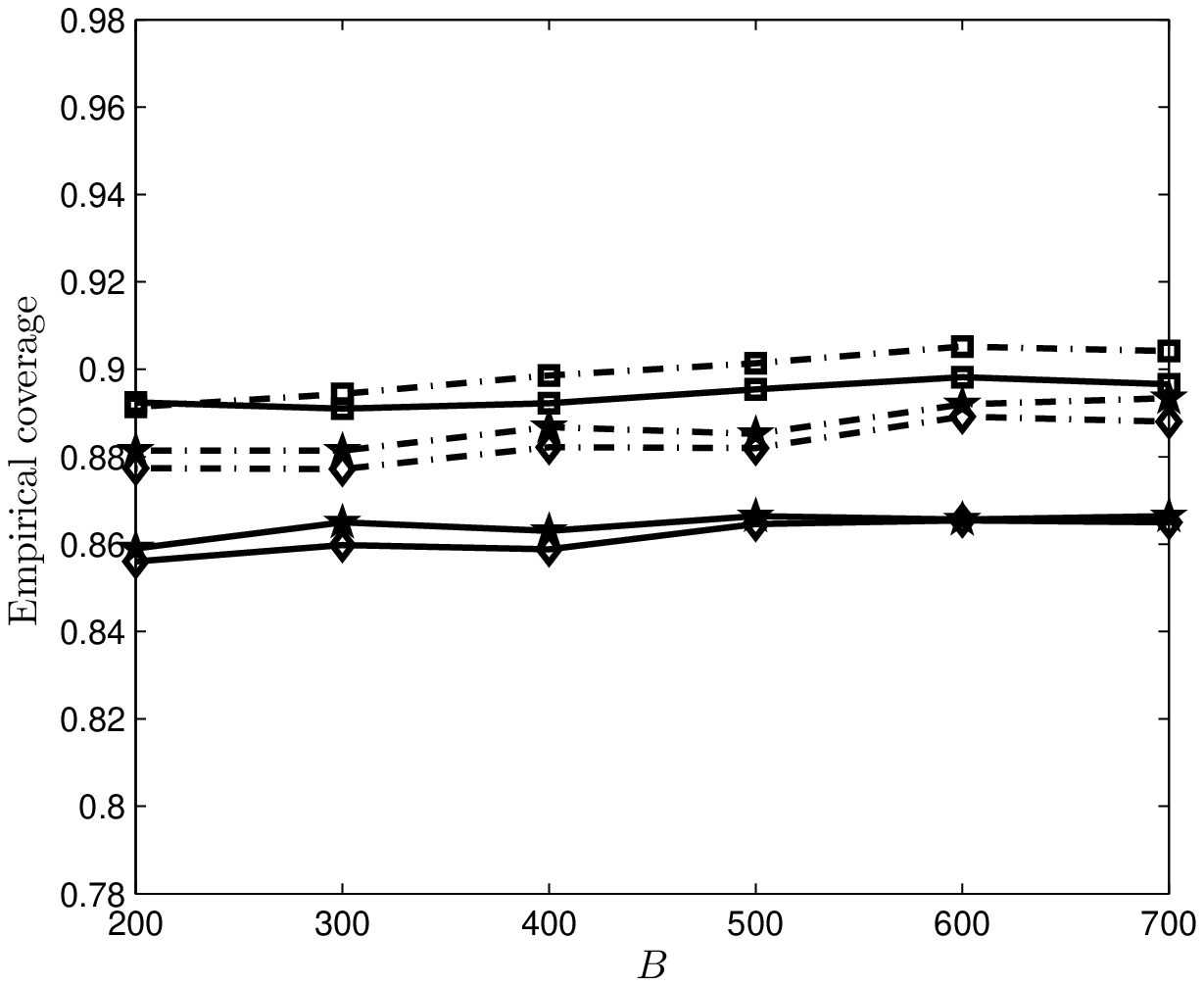}}\\
\subfigure{\includegraphics[scale=0.45]{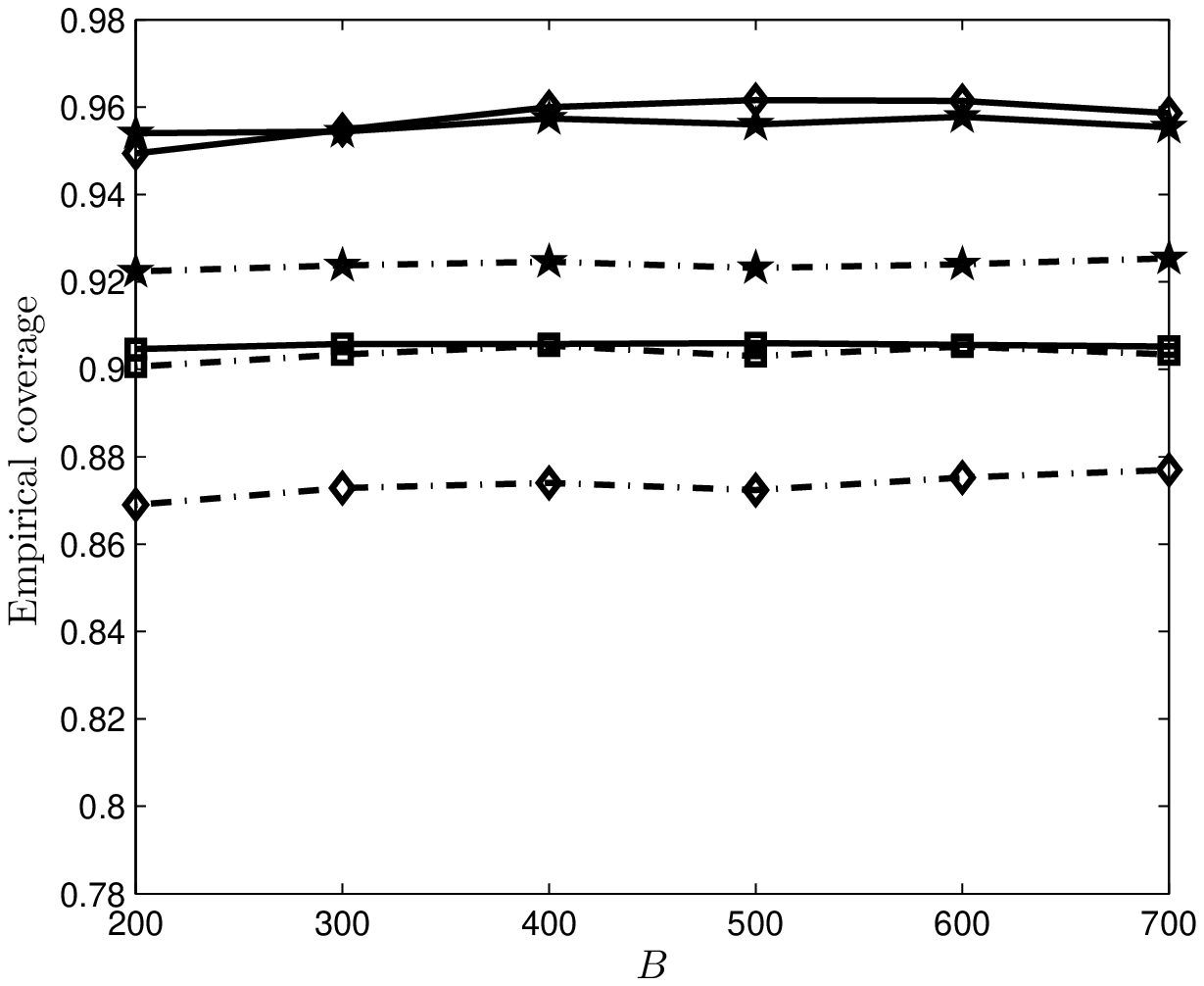}}
\subfigure {\includegraphics[scale=0.45]{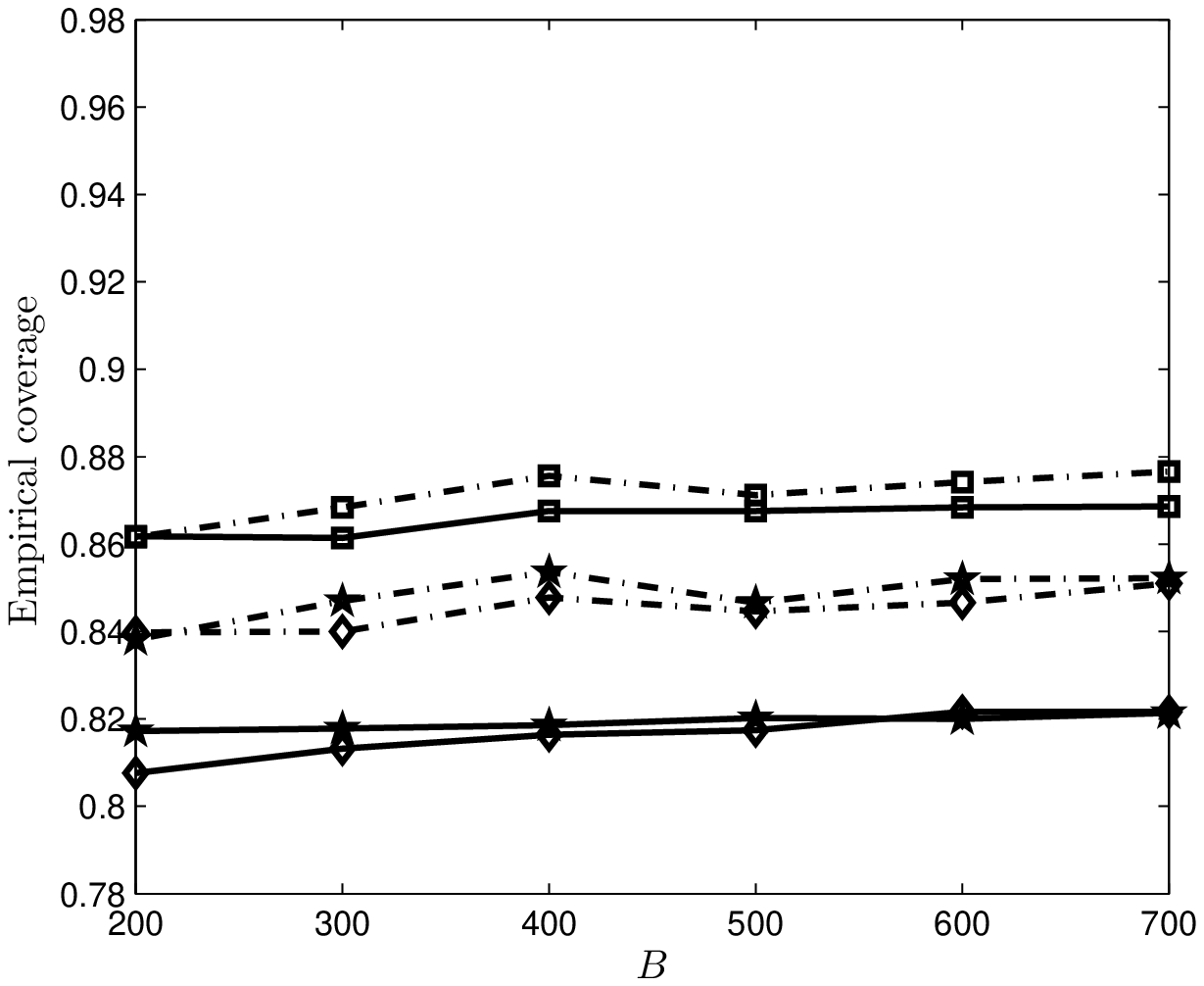}}\\
 \setlength{\abovecaptionskip}{0pt}
\end{center}
\caption{Performance of bootstrap methods for confidence intervals when $n=40$.  First and second rows show results for the exponential distribution, and the log-normal distribution, respectively; left- and right-hand panels show results for one-sided and two-sided equal-tailed confidence intervals, respectively. In each panel the graphs represent single-bootstrap percentile ($-\star-$), single-bootstrap percentile-$t$ ($-\cdot\star\cdot-$), conventional double-bootstrap percentile ($-\square-$), conventional double-bootstrap percentile-$t$ ($-\cdot\square\cdot-$), warp-speed percentile ($-\diamondsuit-$) and warp-speed percentile-$t$ methods ($-\cdot\diamondsuit\cdot-$). }
\end{figure}

\end{document}